\documentclass[11pt, draftcls, onecolumn]{IEEEtran}

\usepackage{latexsym, amssymb,amsfonts,dsfont}
\usepackage{graphicx}\usepackage{subfigure}
\usepackage{tikz}\usepackage{pgf}\usetikzlibrary{arrows,automata}
 \usepackage[latin1]{inputenc}
 \usepackage{verbatim}
\usepackage{amsmath, amsbsy}
\usepackage{amsopn, amstext}
\usepackage{extarrows}
\usepackage{hyperref}
 \usepackage{longtable}
 \usepackage{rotating}
 \usepackage{multirow}
\usepackage{cancel, color} 
\usepackage{cite}
\usepackage{float}
\usepackage{ifpdf} 
 
  \usepackage[justification=centering]{caption}
\newtheorem{thm}{Theorem}[section] 
\newtheorem{lem}[thm]{Lemma}

\newtheorem{rem}[thm]{Remark} 
\newtheorem{ass}{Assumption}

\allowdisplaybreaks[4]

\def\deq{\stackrel{\Delta}{=}}

\title{\LARGE \bf  Asymptotic Properties of   Primal-Dual  Algorithm  for   Distributed  Stochastic Optimization  Over   Random Networks }

\author{Jinlong  Lei, Han-Fu Chen, and Hai-Tao Fang}

\begin{document}
\date{}
\maketitle
\begin{abstract}  This paper studies    a  distributed    stochastic   optimization  problem over random networks with imperfect communications  subject to a global constraint, which is  the intersection of   local constraint sets assigned to   agents. The global cost function is the sum of local cost functions,  each  of which is the expectation of
  a random cost function.  By incorporating   the augmented Lagrange technique  with the projection method,  a   stochastic approximation based  distributed primal-dual algorithm   is proposed to  solve  the problem.
Each agent updates its estimate by using the  local   observations  and the information derived from   neighbors.   
For the constrained problem,     the estimates   are  first shown to be bounded almost surely (a.s.),  and   then    are proved to converge to the optimal solution set  a.s.  Furthermore,   the  asymptotic  normality and efficiency     of the   algorithm are addressed for  the   unconstrained case.  The results  demonstrate the influence of    random  networks,    communication noises, and   gradient errors on the   performance of the algorithm.
    Finally,  numerical simulations demonstrate  the theoretic results. 
\end{abstract}

\begin{IEEEkeywords}
Distributed   stochastic optimization, random networks, primal-dual algorithm, stochastic approximation, asymptotic normality, asymptotic  efficiency. 
\end{IEEEkeywords}

\section{Introduction}  
For recent years,  extensive efforts have been paid to the  distributed   estimation and   optimization  problems   motived by their  wide applications  in sensor networks \cite{WSN1,WSN2},  cognitive networks \cite{WSN3}, multi robots \cite{robot}, as well as in  distributed   learning \cite{Learning}.  This paper  studies  a   distributed     optimization  problem, where $n$ agents connected in a network collectively minimize a convex cost function  $\sum_{i=1}^n f_i(x)$ subject to a convex set constraint $\bigcap_{i=1}^n \Omega_i$. 
 The local cost function of agent  $i$ takes the form $f_i(x)=E[h_i(x,\vartheta_i)],$ where   $\vartheta_i$  is a random variable.  In such a problem,  the   cost function  $f_i(x)$ is difficult to calculate, but   samples of the  cost function $h_i(x,\vartheta_i)$ may  serve as    estimates  for its expectation.  It is assumed that the local constraint set $\Omega_i$ of agent  $i$  is closed and convex, and  the communication relationship  among the  agents is  described by  a random network. Besides,  communications are imperfect since there are noises in the  channels  through which   agents exchange information.

       There exist many  papers   considering the  related problems.  A unconstrained cooperative optimization problem is investigated in  \cite{optimization0}  and  \cite{optimization1}    over  the deterministic and the random switching networks, respectively.     A distributed stochastic subgradient projection algorithm is proposed in \cite{nedic} to solve
       a  constrained optimization problem,   where  all agents are subject to a common convex constraint set, and subgradients of  local cost functions  are  corrupted by  stochastic errors.    
    Effects of stochastic subgradient errors on the convergence of the algorithm  over   deterministic  switching  networks are investigated.  A distributed  asynchronous  algorithm with two diminishing  step sizes 
  is  proposed in  \cite{nediB2} to solve the   distributed constrained  stochastic optimization   problem.
   The estimates  are shown to  converge  to a random point in the optimal set a.s.,  when constraint   sets are  compact, the global constraint set  has  a nonempty interior set, and  cost functions  are non-smooth but  with bounded subgradients.    A distributed algorithm based on dual subgradient averaging is proposed in  \cite{Duchi}, where it is shown how do   the network size  and the spectral gap of the network influence  convergence rates.  Consensus-based distributed primal-dual subgradient methods  are given in\cite{nediB1,zhu}, where \cite{zhu}  solves a problem with the cost function being  the sum of local  cost functions and with global convex inequality constraints known to all agents,  while \cite{nediB1}   solves a problem  with  a coupled global cost function and  with  inequality   constraints. A primal-dual algorithm  with constant step size is  proposed  and its convergence is shown in  \cite{Ling}  for  the deterministic  unconstrained optimization problem  over  undirected  connected  graph with perfect  communications.    Besides, performance of   the continuous time primal-dual  algorithms is also investigated in  \cite{Feijer,Liu}.  Generally  speaking,  the above mentioned  algorithms can be divided into three  categories:     \cite{optimization0,optimization1,nedic, nediB2} belong to   the  primal  domain algorithms, \cite{Duchi} belongs to the dual    domain algorithm,   while  \cite{ nediB1,Ling,zhu,Liu,Feijer} belong to the  primal-dual  domain algorithms.

In this paper,  we propose a  stochastic  approximation based  distributed primal-dual  algorithm     to solve the distributed  constrained stochastic  optimization  problem.   Since  it   is equivalent to a convex optimization problem with a  linear equality constraint  and a convex set constraint, 
 by   incorporating   the augmented Lagrange technique  with the projection method,
  a distributed    primal-dual   algorithm is derived.  The   algorithm is   distributed as  in an iteration each agent updates its   estimate   using   the noisy  observations for gradients of    the  local functions
  and the noisy  observations for  both  primal and   dual variables of  the  neighboring agents.

  Contributions of the  paper are as follows. 1)  
 Stability   and convergence  of the algorithm are proved  for the constrained problem.  The communication  graphs
 are assumed to be independent identically distributed  (i.i.d)  with   the mean graph   being  undirected  and connected.  Communication noises and  gradient errors  are assumed to be martingale difference sequences (mds).   Convex sets are required to have smooth boundaries  with the global constraint  having at least one relative interior, and  gradient  functions are required to be Lipschitz continuous.   Then with diminishing step-size,   the estimates   are shown to be bounded  a.s. by  using the  convergence theorem for martingales, and  to  converge to the optimal solution set a.s.  by  use of  the results for constrained stochastic approximation \cite{Kushner}.       Compared with    \cite{nediB2},  here   gradients  are only  required to be Lipschitz continuous  without   boundedness assumption,  constraint sets are not assumed to be compact, and   the global constraint is only  required to have at least one  relative interior point.    Different from   \cite{ nediB1,Ling,zhu},     the stochastic optimization  problem   is investigated  over random networks  with  imperfect communications.  2)  Asymptotic properties are considered for the unconstrained problem.     Through dimensionality  reduction,  asymptotic normality and efficiency   of the algorithm  are established.   In comparison with  \cite{Duchi}, we  have shown  the exact influence of    random networks, imperfect communications,  gradient errors   and the structure of cost functions on  the  rate of convergence.

       The organization  of the paper is as follows.   In Section \ref{sec:PS},
 some preliminary information about   graph theory and convex analysis is provided and the problem is formulated.
 In Section \ref{sec:Algorithm},    the basic  assumptions are introduced  and  a  stochastic approximation based distributed  primal-dual    algorithm is designed. Convergence  for the constrained problem is  established in Section \ref{sec:Convergence}, 
  while asymptotic properties for the unconstrained problem are given  in Section  \ref{sec:Normality}.
 Numerical  simulations are demonstrated  in Section \ref{sec:Simulation} with  some concluding remarks given   in Section \ref{sec:Conclusion}.

\section{Preliminaries and Problem Statement} \label{sec:PS}
We first introduce  some preliminary  results   about graph theory, convex functions and convex sets,
 then formulate the distributed    optimization problem.  
 \subsection{ Graph Theory } 
Consider a network of  $n$  agents.  The communication relationship  among 
  agents is described by a digraph   $\mathcal{G } =(\mathcal{V }, \mathcal{E }_{\mathcal{G }},\mathcal{A}_{\mathcal{G }}\}$, where  $\mathcal{V }=\{ 1,\cdots,n\}$ is the node set 
   with  node  $i$ representing  agent  $i$;
 $\mathcal{E }_{\mathcal{G }}  \subset  \mathcal{V } \times \mathcal{V } $ is the edge set,
  and  $(j,i)\in\mathcal{E }_{\mathcal{G }}$ if and only if  agent  $i$ can get information from agent  $j$;
 $\mathcal{A}_{\mathcal{G }} =[a_{ij} ]\in  \mathbb{R}^{n\times n} $ is the adjacency matrix of $\mathcal{G } $,
where  $a_{ij} >0$  if  $(j, i)\in \mathcal{E }_{\mathcal{G }} $, and  $a_{ij}=0$, otherwise. 
 Here, we assume the self-edge $(i,i)$ is not  allowed, i.e., $a_{ii}=0~\forall i \in \mathcal{V } .  $
 The Laplacian matrix of   graph   $\mathcal{G}$ is defined as
 $\mathcal{L}_{\mathcal{G}}= \mathcal{D}_{\mathcal{G}}-\mathcal{A}_{\mathcal{G}}$ with $\mathcal{D}_{\mathcal{G}}=diag\{ \sum_{j=1}^n a_{1j}, \cdots, \sum_{j=1}^n a_{nj})$,  where and hereafter  
$diag\{D_1,\cdots, D_n\}$ denotes  the  block diagonal matrix  with   main diagonal blocks being  square matrices  $D_i,~i=1,\cdots,n,$ and with the off-diagonal blocks being  zero matrices. 

  For a   bidirectional graph $\mathcal{G } $,
  $ (i,j)\in \mathcal{E }_{\mathcal{G }}   $ if and only if $ (j,i)\in \mathcal{E }_{\mathcal{G }}$.
The graph  $\mathcal{G } $ is   undirected    if   $\mathcal{A}_{\mathcal{G }}$ is symmetric. 
  The undirected  graph  $\mathcal{G } $  is connected  if for any pair  $i,j\in \mathcal{V}$, there exists  a sequence of nodes
 $i_1, \cdots, i_{p } \in \mathcal{V}$ such  that $ (i,i_1)\in \mathcal{E }_{\mathcal{G }} , $ $ (i_1,i_2)\in \mathcal{E }_{\mathcal{G }} $,
  $\cdots$, $(i_{p},j) \in \mathcal{E }_{\mathcal{G }} $.           For   matrix  $A=[a_{ij}] \in \mathbb{R}^{n \times n}$ with $a_{i j}\geq 0~ \forall i,j=1,\cdots,n$, denote by
    $\mathcal{G}_{A}= \{\mathcal{V},\mathcal{E}_{\mathcal{G}_{A}},  \mathcal{A}_{\mathcal{G}_{A}}\}$  the digraph generated by $A$, where  $\mathcal{V }=\{ 1,\cdots,n\}$,  $ \mathcal{A}_{\mathcal{G}_{A}}= A$, and 
   $(j,i)\in\mathcal{E}_{\mathcal{G}_{A}}$ if   $a_{ij} >0 $.

 The following lemma presents  some  properties of  the Laplacian matrix $\mathcal{L}  $  corresponding to an
  undirected   graph $\mathcal{G}  $.

\begin{lem}  \label{lem1} \cite{graph}
The Laplacian  matrix   $\mathcal{ L}  $   of  an  undirected graph $\mathcal{G} $ has the following properties:
  
  i)   $\mathcal{L}$   is symmetric and positive semi-definite;
  
 ii)   $\mathcal{ L}$  has a simple zero eigenvalue and the corresponding eigenvector is  $\mathbf{1}$,  
 and all the other eigenvalues are positive if and only if  $\mathcal{G} $ is connected,
 where  $\mathbf{1}$  denotes the vector  with all entries equal to 1.

\end{lem}

\subsection{Gradient, Projection  Operator  and  Normal Cone}
For  a given function $f: \mathbb{R}^m \rightarrow [-\infty, \infty],$
  denote its  domain  as   $ \textrm{dom} (f) \triangleq(x \in \mathbb{R}^m:  f(x) < \infty\}.$
Let $f(\cdot)$ be a convex function, and let $x \in \textrm{dom} (f)$.
For  a smooth (differentiable) function $f(\cdot)$,  denote  by $\nabla f(x)$ and by  $\nabla^2 f(x)$ the gradient and Hessian of  $f(\cdot)$ at point $x$, respectively. 
Then  \begin{equation} \label{gradient}
f(y) \geq f(x) + \nabla f(x)^T (y-x) ~~  \forall y \in \textrm{dom} (f),
\end{equation}
 where $x^T$ denotes the transpose of $x.$

 For a  nonempty closed convex   set $\Omega \subset \mathbb{R}^m$ and a point $x \in \mathbb{R}^m$, we call the point in $\Omega$ that is closest  to $x$ the projection of $x$ on $\Omega$ and   denote it by $P_{\Omega} (x)$.  $P_{\Omega} (x)$ contains only one element for any $x \in  \mathbb{R}^m,$ 
 and it satisfies the following non-expansive property  \cite[Theorem 2.13]{opt}
   \begin{equation}\label{pro} 
\| P_{\Omega} (x)-P_{\Omega} (y) \| \leq \| x-y\| ~~ \forall x, y\in\mathbb{R}^m.
\end{equation} 
Consider  a convex closed set $\Omega \subset \mathbb{R}^m$ and a point $x \in  \Omega$. Define
   the normal cone  to $\Omega$ at $x$ as   $N_{\Omega}(x) \triangleq(v \in \mathbb{R}^m: \langle v, y-x \rangle \leq 0  ~~\forall y \in \Omega\}$.  It is shown that  \cite[Lemma 2.38]{opt}  
 \begin{equation}\label{normalcone}
N_{\Omega}(x) =\{ v \in \mathbb{R}^m: P_{\Omega}(x+v)=x)~~\forall x \in  \Omega.
\end{equation}

   A set  $C$ is affine if it contains the lines that pass through  any
 pairs of points $x,y \in C$ with $x \ne y$. Let   $\Omega \subset \mathbb{R}^m$ be a nonempty convex set.
  We say that $x \in \mathbb{R}^m$ is   a   relative interior point of   $\Omega $ if 
 $x \in \Omega$ and  there exists an open sphere $S$ centered at $x$  such  that
$S \cap \textrm{aff}(\Omega) \subset \Omega,$
where   $\textrm{aff}(\Omega)$  is the intersection of all affine sets containing $\Omega$. 
A  pair of vectors  $ x^{*} \in  \Omega  $ and $  z^{*} \in  \Psi$      is  called a saddle point
 of the function $ \Phi(x,  z)$ in $ \Omega \times \Psi$ if  \begin{equation}
\Phi(x^{*}, z) \leq \Phi(x^{*}, z^{*})  \leq \Phi(x , z^{*}) ~~ \forall  x \in \Omega,~ ~ \forall  z \in \Psi. \nonumber
\end{equation}
    These  definitions can be found in \cite{Bertsekas}.

\subsection{Problem Statement}  
 Consider a network of $n$ agents. The objective of the network is to  solve the  following constrained optimization   problem 
\begin{equation}\label{problem1}
\begin{split}
&  \textrm{minimize} ~~ f(x)=\sum_{i=1}^n  f_i(x) ,\\
&  \textrm{subject to} ~~ x \in  \Omega_o= \bigcap_{i=1}^n \Omega_i,
\end{split}
\end{equation}
where $f_i(x): \mathbb{R}^m  \rightarrow \mathbb{R}$ is  the local
cost function of agent $i$, and $\Omega_i \subset \mathbb{R}^m$
  is the local constraint set of agent  $i$. Assume that $f_i(\cdot)$ is a smooth convex function on  $\Omega_i$,
    and  $\Omega_i$ is  a closed convex set   only  known to  agent  $i$.  
   Assume   there  exists at least one finite solution $x^*$ to the problem \eqref{problem1}.
For   the problem \eqref{problem1},  denote by   $f^{*}= \min_{ x \in     \Omega_o} f(x)$ the optimal value,
and  by $  \Omega_{o}^{*}=\{  x \in   \Omega_o: f(x)=f^{*}) $  the optimal solution set.

Further, assume that for each $i \in \mathcal{V}$,  the values of $f_i(\cdot)$ and $\nabla f_i(\cdot)$ are observed with  noises.  For example,  $f_i(x)=E[h_i(x,\vartheta_i)],$ where $h_i:\mathbb{R}^m \times \Theta_i \rightarrow \mathbb{R}$ with    $\vartheta_i$  being a random variable defined on $\Theta_i$, and  the expectation  $E[\cdot]$ is taken with respect to $\vartheta_i.$ 
In  this case, one may only observe    $\nabla h_i(x_i, \vartheta_i)$  for some given  samples of $\vartheta_i$, 
while the  exact gradient $\nabla f_i(x_i)$  is difficult to calculate.  

Let  the communication relationship among agents   at time $k$  be  described by  a  directed  graph 
  $\mathcal{G}_k =\{  \mathcal{V}, \mathcal{E}_{\mathcal{G}_k}, \mathcal{A}_{\mathcal{G}_k}\}$,  where
    $\mathcal{V }=\{ 1,\cdots,n\}$ is the node set,   $\mathcal{E}_{\mathcal{G}_k}$ is the edge set, and  $\mathcal{A}_{\mathcal{G}_k} =[ a_{ij,k} ]_{i,j=1}^n $ is the adjacency matrix. 
 Denote by   $\mathcal{L}_k=[ l_{ij,k} ]_{i,j=1}^n$  the   Laplacian  matrix  of digraph $ \mathcal{G}_k$. 
  Denote by $\mathcal{N}_{i,k}=\{j \in \mathcal{V}:(j,i) \in  \mathcal{E}_{\mathcal{G}_k}\}$ the neighbors  of 
  agent  $i$ at time $k.$ Besides,   neighboring agents  exchange  information through channels which may contain noises. 
The noises  may be introduced  by  quantization errors \cite{Xie,Rabbat},  or actively introduced   to achieve differential privacy  \cite{DCU1}.

   \section{Primal-Dual Algorithm  } \label{sec:Algorithm}
  We now propose a  distributed  primal-dual algorithm to solve the distributed stochastic optimization  problem, and list   some conditions  and   preliminary lemmas  to be used in the sequel. 
\subsection{Algorithm Design }

Denote  by $x_{i,k} \in \mathbb{R}^m$ the estimate for the  optimal  solution to problem \eqref{problem1} given by agent  $i$ at time $k$,  and by $\lambda_{i,k} \in \mathbb{R}^m$    the auxiliary variable  of  agent  $i$. Hereafter, 
we call $x_{i,k}$  and $\lambda_{i,k}$  the   primal and dual variables for  agent  $i$ at time $k$.  Agents exchange  information in  the following way: 
  if $(j,i) \in  \mathcal{E}_{\mathcal{G}_k}$,  then agent  $i$    gets the noisy observations  $\{x_{ij,k}, \lambda_{ij,k}\}$ of   $\{x_{j,k}, \lambda_{j,k}\}$    given as  follows:
   \begin{equation}
\label{ }
\begin{split}
& x_{ij,k} =  x_{j,k}+ \omega_{ij,k} \textrm{~ if  }(j,i) \in  \mathcal{E}_{\mathcal{G}_k},   \\& 
  \lambda_{ij,k} =  \lambda_{j,k}+ \zeta_{ij,k} \textrm{~ if  }(j,i) \in  \mathcal{E}_{\mathcal{G}_k}, 
 \end{split}
\end{equation}
  where $\omega_{ij,k}$ and $\zeta_{ij,k}$  denote  the communication noises.
  
   The    sequences $\{  x_{i,k}\}$  and    $\{ \lambda_{i,k}\}$ are updated  as follows:
 \begin{equation}\label{algorithm1}
\begin{split}
 &  x_{i , k+1 }= P_{\Omega_i} \big( x_{i,k}- \gamma_k g_{i,k}  
  - \gamma_k \sum_{j \in \mathcal{N}_{i,k}}   a_{ij,k}(\lambda_{i,k} -\lambda_{ij,k})   - \gamma_k \sum_{j \in \mathcal{N}_{i,k}}   a_{ij,k}(x_{i,k} -x_{ij,k})  \big), \\
& \lambda_{i ,k+1}=\lambda_{i,k}+\gamma_k\sum_{j\in \mathcal{N}_{i,k}}  a_{ij,k} (x_{i,k}-x_{ij ,k}),
\end{split}
\end{equation}  
where  $\gamma_k$ is the step size and $g_{i,k}$ denotes the noisy observation of $ \nabla  f_i(x_{i,k}) $:
\begin{equation}
\label{ }
g_{i,k}=  \nabla  f_i(x_{i,k}) +  v_{i,k},
\end{equation}
 where $v_{i,k}$ is the observation  noise.
   Note that the algorithm \eqref{algorithm1} is distributed as  in an  iteration each agent updates its  local estimates only using  the  local  gradient  observations  and   the noisy observations for primal and dual variables of its neighbors.   

Set  $ X_k \triangleq  col\{  x_{1,k} , \cdots, x_{n,k}\}, $ $ \Lambda_k \triangleq  col\{ \lambda_{1 ,k} \cdots, \lambda_{n ,k}\},$
and $   \nabla \widetilde{f} (X_k) \triangleq  col\{   \nabla f_1(x_{1,k})   ,\cdots,   \nabla f_n(x_{n, k}) \},$ 
where  by $ col\{x_1,\cdots, x_n\} $  we mean $ (x_1^T, \cdots, x_n^T)^T$.   Define 
 $v_k \triangleq  col\{   v_{1,k} ,\cdots, v_{n,k}\}$,   $\omega_k \triangleq  col\{  \omega_{1,k} ,\cdots, \omega_{n,k}\}$
with $\omega_{i,k} \triangleq  \sum_{j=1}^n a_{ij,k}\omega_{ij,k}$,
and $\zeta_k \triangleq  col\{ \zeta_{1,k},\cdots, \zeta_{n,k}\}$ with $\zeta_{i,k}\triangleq  \sum_{j=1}^n a_{ij,k}\zeta_{ij,k}$. 
Then the  algorithm \eqref{algorithm1} can be rewritten in the  compact form as follows
 \begin{equation}\label{compact0}
\begin{split}
& X_{k+1}= P_{\Omega} \big( X_k -  \gamma_k  \nabla \widetilde{f}(X_{k})  -\gamma_k( \mathcal{ L} _k \otimes\mathbf{I} _m)  \big( \Lambda_k +  X_k \big)   + \gamma_k  \big(  \zeta_k +\omega_k-v_k \big) \big), \\
&\Lambda_{ k+1}=\Lambda_k+\gamma_k   ( \mathcal{L}_k \otimes\mathbf{I} _m)  X_k- \gamma_k \omega_k  ,
\end{split}
\end{equation}
where 
$\Omega=\prod_{i=1}^n \Omega_i$ denotes the Cartesian product, the symbol 
$\otimes $ denotes  the Kronecker product,  and  $\mathbf{I}_m$  denotes the identify matrix of size $m$.

\subsection{Assumptions}

We impose the  following        assumptions on the constraint sets and on the   cost functions.
\begin{ass}\label{ass-set}

a) $ \Omega_o$  has at least one relative   interior point.

 b) There exists a constant $L_f>0$ such that  for any $i \in \mathcal{V}$
 \begin{equation}\label{Lipschitz}
  \| \ \nabla  f_i(x)  - \nabla f_i(y)\| \leq  L_f \| x-y \|~~ \forall  x,y \in \Omega_i.
\end{equation}

 c) For any $i \in \mathcal{V}, $ the set $\Omega_i$ is determined by $p_i$ inequalities:
$$ \Omega_i=\{ x\in \mathbb{R}^m: q_{ij} (x) \leq 0, ~ \forall  j=1,\cdots, p_i\},$$
where   $q_{ij}(\cdot),~j=1,\cdots, p_i$ are continuously differentiable real-valued functions on
$\mathbb{R}^m$.
Moreover,  $\{ \nabla q_{ij}(x), ~ j \in A_i(x)\}$ are linearly independent,  where $A_i(x)=\{ j: q_{ij}(x)=0\}$.
\end{ass}  
\begin{rem}
The existence of  the relative   interior point   will be used   to  guarantee  that   the  primal and dual problems defined in  Section \ref{sec:Algorithm}. C have the same optimal solution.  The globally Lipschitz  condition is used to guarantee the boundedness of the estimates. Assumption \ref{ass-set}-c   indicates  that  
all  local constraint sets   have   smooth boundaries.  In fact,    Assumption \ref{ass-set}-c corresponds to   A4.3.2 in  \cite{Kushner}  but without compactness  requirement.   \end{rem}

The   following  conditions are imposed  on the communication graphs and on the adjacency matrices.

\begin{ass} \label{ass-graph}(Mean graph is connected and undirected)

 a)   $(\mathcal{A}_{\mathcal{G}_k} )_{ k \geq 0}$ 
is an i.i.d   sequence with expectation   denoted by  $\mathcal{\bar{A}} = E [  \mathcal{A}_{\mathcal{G}_k} ]$.

 b)   The   graph    $\mathcal{G}_{ \mathcal{\bar{A}}}$ generated by $\mathcal{\bar{A}} $ is  undirected and connected.

 c)  There  exists a constant $\eta>0$  such that 
$$ E [  a_{ij,k} ^2] =\sigma_{ij}\leq  \eta^2 ~~ \forall i,j \in \mathcal{V}.$$ 

  d) $\mathcal{L}_k$ is independent of $\mathcal{F}_{k-1}$,  where 
\begin{equation}\label{def-algebra}
\mathcal{F}_{k}  = \sigma\{X_0, \Lambda_0, \omega_{ij,t}, \zeta_{ij,t}, v_{i,t},
\mathcal{L}_t, 0 \leq t \leq k , 1\leq i,j \leq n\}.
\end{equation}   
 \end{ass}

 \begin{rem}
 Note that   Assumption \ref{ass-graph} does not  require the random   graph  
 at  any instance    be  undirected or  strongly  connected.  It only requires   the mean graph  be  
 undirected and connected.  The gossip-based communication protocol \cite{gossip}
 and   the broadcast-based communication  \cite{broadcast} 
 both  satisfy   Assumption \ref{ass-graph}  when  the underling graph is  bidirectional and strongly connected. 
\end{rem}

Set 
\begin{equation}\label{def-algebras}
\mathcal{F}_{k}'= \sigma\{ \mathcal{L}_{k+1 } , \mathcal{F}_{k}\}.
\end{equation}

Note that the adjacency matrix  $\mathcal{A}_{\mathcal{G}_k} $  is  uniquely defined by $ \mathcal{L}_{k } $ with
$a_{ij,k}=-l_{ij,k}~\forall i\neq j $ and $ a_{ii,k}=0$.  Thus, the covariance of  $ \mathcal{L}_{k } $ 
is finite  by Assumption \ref{ass-graph}-c,    $\mathcal{A}_{\mathcal{G}_k} $  is independent of $\mathcal{F}_{k-1}$ by Assumption \ref{ass-graph}-d, and  $ \mathcal{L}_{k } $  is adapted to $\mathcal{F}_{k-1}' $ by its definition \eqref{def-algebras}.

The following conditions are imposed   on  the communication  noises and    gradient errors.

\begin{ass}\label{ass-noise2} 

 a)  For any $i,j\in \mathcal{V},$  $\{\omega_{ij,k} ,\mathcal{F}_{k}'\}$ is an mds  with
 \begin{align}
& E[\omega_{ij,k} |\mathcal{F}_{k-1}']=\mathbf{0},~~ E[\| \omega_{ij,k} \|^{2 } |\mathcal{F}_{k-1}'] \leq \mu^2,\nonumber
 \end{align} 
 and
 \begin{align} 
&   E[\omega_{ij,k}  \omega_{ij,k} ^T |\mathcal{F}_{k-1}'] \triangleq R_{\omega,ij} .\label{covW} \end{align} 
   
b)  For any $i,j\in \mathcal{V},$ $\{\zeta_{ij,k} ,\mathcal{F}_{k}'\}$ is an mds with
$$ E[\zeta_{ij,k} |\mathcal{F}_{k-1}']=\mathbf{0}~~E[\| \zeta_{ij,k} \|^{2 } |\mathcal{F}_{k-1}'] \leq \mu^2,$$
and  
 \begin{align} &   E[ \zeta_{ij,k}  \zeta_{ij,k} ^T  |\mathcal{F}_{k-1}'] \triangleq R_{\zeta,ij}. \label{covZ} 
\end{align}      

c) For any $i\in \mathcal{V}$,  $\{v_{i,k},\mathcal{F}_{k}'\}$   is an mds  with 
 \begin{align}
& E[v_{i,k} | \mathcal{F}_{k-1}' ]=\mathbf{0}, ~~E[ \| v_{i,k} \|^{2 } | \mathcal{F}_{k-1}'  ] \leq  c_v(1+\|x_{i,k}\|^{2 }), \label{gradientnoise}
 \\
& \lim_{k \rightarrow \infty}  E[v_{i,k} v_{i,k} ^T  |\mathcal{F}_{k-1}'] \triangleq R_{v,i}. \label{covV} \end{align}
\end{ass}

In Section \ref{sec:Convergence}, \eqref{covW}, \eqref{covZ} and \eqref{covV} are not needed. 
The simplified version of  Assumption \ref{ass-noise2}  with  \eqref{covW}, \eqref{covZ} and \eqref{covV} removed will be called Assumption \ref{ass-noise}.  

\begin{ass}  \label{ass-noise} 

a) 
  Assumption \ref{ass-noise2}-a  with  \eqref{covW}  removed. 
  
  b) 
  Assumption \ref{ass-noise2}-b  with \eqref{covZ}   removed. 
  
  c)
  Assumption \ref{ass-noise2}-c  with   \eqref{covV} removed. 
\end{ass}

 \begin{rem}
 The communication noises  introduced by   the probabilistic quantization \cite{Xie,Rabbat} is shown to be an i.i.d sequence with bounded second moments,  and hence 
 satisfy  Assumption \ref{ass-noise}-a and   \ref{ass-noise}-b.  
 Assumption \ref{ass-noise2}-c holds true in many cases, for example, in the   quadratic distributed  stochastic optimization  problem \eqref{filter1}  discussed  in Section \ref{sec:Simulation}. 
 
  \end{rem}  
    
By Assumption \ref{ass-graph}-a,  $\{\mathcal{L}_k \}_{ k \geq 0}$  is an i.i.d sequence. Set  $\mathcal{\bar{L}}\triangleq  E [\mathcal{L}_k ]$.   Then  $\mathcal{\bar{L}} $ is the  Laplacian matrix of  the undirected connected graph  $\mathcal{G}_{ \mathcal{\bar{A}}}$. 
Define  
\begin{align}
& e_{1,k}\triangleq \big( ( \mathcal{\bar{L}}-\mathcal{ L}_k  ) \otimes\mathbf{I} _m \big)   ( \Lambda_k +  X_k ), \label{def-e1}
 \\&e_{2,k}\triangleq \zeta_k +\omega_k-v_k,\label{def-e2}
  \\&e_{3,k}\triangleq  \big(( \mathcal{L}_k -\mathcal{\bar{L}}) \otimes\mathbf{I} _m \big)  X_k-   \omega_k. \label{def-e3}
  \end{align}
   Then \eqref{compact0} can be rewritten as:
 \begin{equation}\label{compact}
\begin{split}
& X_{k+1}= P_{\Omega} \big( X_k -  \gamma_k     \nabla \widetilde{f}(X_{k})   -\gamma_k ( \mathcal{\bar{L}}  \otimes\mathbf{I} _m)  \big( \Lambda_k +  X_k \big)   + \gamma_k  \big(e_{1,k}+e_{2,k}\big) \big), \\
&\Lambda_{k+1}=\Lambda_k+\gamma_k   ( \mathcal{\bar{L}}  \otimes\mathbf{I} _m)  X_k+  \gamma_k  e_{3,k}  . 
\end{split}
\end{equation}

We impose the following condition on the step size $\{\gamma_k\}$.

\begin{ass}\label{ass-stepsize} 
$$\gamma_k >0, ~ \sum_{k=1}^{\infty}  \gamma_k =\infty,\textrm{  and }\sum_{k=1}^{\infty}  \gamma_k^2 < \infty .$$
\end{ass}

  \subsection{Preliminary Lemmas}
  We now give some preliminary results about   the formulated distributed   optimization  problem. 

\begin{lem} \label{lem2} 
The problem \eqref{problem1} is equivalent to the following constrained optimization   problem 
\begin{equation}\label{problem2}
\begin{split}
& \textrm{minimize}  ~~\widetilde{f}(X)\deq \sum_{i=1}^n  f_i(x_i) ,\\
& \textrm{subject to} ~~ ( \mathcal{\bar{L}} \otimes\mathbf{I} _m ) X=\mathbf{0}, ~~  X \in \Omega ,
\end{split}
\end{equation}
where $X= col\{x_1,\cdots, x_n\} $.\end{lem}

The result can be easily derived since  $( \mathcal{\bar{L}} \otimes\mathbf{I} _m ) X=\mathbf{0}$
 if and only if $x_i=x_j~\forall i,j\in \mathcal{V}.$

Define $ \Phi(X,\Lambda)\triangleq \widetilde{f}(X)+\Lambda^T ( \mathcal{\bar{L}} \otimes\mathbf{I} _m ) X$  as the Lagrange function, where  $\Lambda \in  \mathbb{R}^{mn}$ is the Lagrange  multiplier.   Then the  problem  \eqref{problem2} can   be rewritten as $ \inf\limits_{X \in \Omega} \sup\limits_{\Lambda \in \mathbb{R}^{ mn} } \Phi(X, \Lambda)$, while  the  dual problem is defined as follows
\begin{equation}\label{pdual}
\sup\limits_{\Lambda \in \mathbb{R}^{ mn} } \inf\limits_{X \in \Omega}  \Phi(X, \Lambda).
\end{equation}

\begin{lem}   \label{lem3} Assume Assumption \ref{ass-set}-a and  Assumption \ref{ass-graph}-b  hold.    Then $ \Phi(X, \Lambda)$ has at least one saddle point in $ \Omega \times \mathbb{R}^{mn}$.
 A  pair  $(X^{*}, \Lambda^{*}) \in  \Omega \times \mathbb{R}^{mn}$  is
the primal-dual  solution   to the problems  \eqref{problem2} and \eqref{pdual}
  if and only if    $(X^{*}, \Lambda^{*}) $ is  a saddle point
 of $ \Phi(X, \Lambda)$ in $ \Omega \times \mathbb{R}^{mn}$.
\end{lem}   
{\bf Proof}:
   Assumption \ref{ass-set}-a  implies that there exists  a  relative interior $\bar{X} $ of  set
$\Omega$ such that $(\mathcal{L}  \otimes \mathbf{I}_m) \bar{X}=0$.  Since $f^{*}$ is finite,    by    \cite[Proposition   5.3.3]{Bertsekas} we know that  
 \begin{equation}\label{minmax}
\inf_{X \in \Omega} \sup_{\Lambda \in \mathbb{R}^{ mn} } \Phi(X, \Lambda)=\sup_{\Lambda \in \mathbb{R}^{ mn} } \inf_{X \in \Omega}  \Phi(X, \Lambda),
\end{equation}and there  exists  at least one dual optimal solution.  

Since the  minimax equality \eqref{minmax} holds,   by   \cite[Proposition 3.4.1]{Bertsekas}  
 $X^{*}$ is  the  primal optimal   solution and $\Lambda^{*}$  is the dual optimal   solution if and only if
 $(X^{*}, \Lambda^{*})$ is  a saddle point  of $ \Phi(X,\Lambda)$ on  $\Omega \times \mathbb{R}^{mn}$.
Since there exists at least one  primal and dual optimal solution pair, we conclude that $ \Phi(X, \Lambda)$ has at least one saddle point in $ \Omega \times \mathbb{R}^{mn}$.
  This completes the proof. \hfill $\blacksquare$

\section{Convergence Theorems} \label{sec:Convergence}
In this section, we analyze   stability  and convergence of    the  algorithm \eqref{algorithm1}.  
For  notational simplicity, we     assume   $m=1$  in this section. This  does not influence the convergence    analysis for the general case $m \geq 1$. 
\subsection{Stability Analysis}
     \begin{thm}\label{thm2} (Stability)  Let $\{ x_{i,k}\}$ and  $\{ \lambda_{i,k}\}$   be produced by the algorithm \eqref{algorithm1} with any   initial values $x_{i,0},~ \lambda_{i,0}$. 
 Let Assumptions \ref{ass-set}-a, \ref{ass-set}-b,   \ref{ass-graph},   \ref{ass-noise},     and  \ref{ass-stepsize}  hold. Then  $  \| X_k -X^{*}\|^2+ \| \Lambda_k-\Lambda^{*} \|^2 $   converges a.s.,
 where $(X^{*}, \Lambda^{*})$   is a saddle point  of $ \Phi(X, \Lambda)$ in $ \Omega \times \mathbb{R}^{mn}$. 
  \end{thm}
 
This theorem establishes that the sequences $\{ X_k\}$ and $\{\Lambda_k\}$ are bounded a.s.,
and the distance  between the pair $(X_k, \Lambda_k)$ and  the saddle point $(X^{*}, \Lambda^{*})$  converges a.s. 
 Before proving the theorem, we first give some preliminary  lemmas.
The following lemma  establishes    properties of noise sequences $\{e_{1,k})$,  $\{e_{2,k}\}$ and $\{e_{3,k}\}$ 
 defined in   \eqref{def-e1},  \eqref{def-e2}, and  \eqref{def-e3}, respectively.

\begin{lem}  \label{lemma1}
Let    Assumptions \ref{ass-graph}-a,   \ref{ass-graph}-c,   \ref{ass-graph}-d,  and 
  \ref{ass-noise}  hold. Then   the following  assertions   take place  a.s. \begin{align}
&E[ e_{1,k} | \mathcal{F}_{k-1}]=\mathbf{0},~~E[ \| e_{1,k}\|^2 | \mathcal{F}_{k-1}]  
   \leq C_{01} \big\| \Lambda_k+  X_k  \big\|^2,\label{noise1}  \\&
E[ e_{2,k} | \mathcal{F}_{k-1}]=\mathbf{0}, ~~E[ \| e_{2,k} \|^2 | \mathcal{F}_{k-1} ]   =C_{02} +3c_v \|X_k\|^2,\label{noise2}  \\&
E[ e_{3,k} | \mathcal{F}_{k-1}]=\mathbf{0}, ~~E [ \| e_{3,k} \|^2 | \mathcal{F}_{k-1}]\leq  C_{01} \| X_k \|^2 + C_{03}, \label{noise3} 
\end{align}
where  $C_{01}=E [\|  \mathcal{L}_k -\mathcal{\bar{L}}\|^2]$,  $C_{02}=3c_v n+  6n^3 \mu^2 \eta^2$,
and $C_{03}=n^3\mu^2 \eta^2$.
\end{lem}
\textbf{Proof:}  By Assumption \ref{ass-graph}-d  we have 
\begin{equation}\label{conditional-Laplacian}
E \big[    \mathcal{\bar{L}}-\mathcal{ L} _k 
 | \mathcal{F}_{k-1}  \big]  = \mathcal{\bar{L}}- E [ \mathcal{ L}_k ]= \mathbf{0}.
\end{equation}

 Since  $X_k $ and $\Lambda_k$  are adapted to $\mathcal{F}_{k-1}  $
  by  \eqref{algorithm1}   \eqref{def-algebra}, from  \eqref{def-e1} \eqref{conditional-Laplacian} it follows that 
\begin{align}
&E[ e_{1,k} | \mathcal{F}_{k-1}]=    E \big[   \mathcal{\bar{L}}-\mathcal{ L} _k  
 | \mathcal{F}_{k-1}  \big]    ( \Lambda_k+  X_k )=\mathbf{0}, \nonumber \end{align}
and  
\begin{equation}
\begin{split}
E[ \| e_{1,k}\|^2 | \mathcal{F}_{k-1}] & \leq    E \big[   \| \mathcal{\bar{L}}-\mathcal{ L} _k   \|^2
 | \mathcal{F}_{k-1}  \big]  \cdot  \| \Lambda_k+  X_k   \|^2\\&
   =E [\|  \mathcal{L}_k -\mathcal{\bar{L}}\|^2]  \cdot  \| \Lambda_k+  X_k   \|^2.  \nonumber
   \end{split}\end{equation}
 Therefore,  \eqref{noise1} holds.

Since $a_{ij,k}$ is adapted to  $\mathcal{F}_{k-1}'$ by \eqref{def-algebras}, 
    from Assumption \ref{ass-noise}-a it follows   that  for any $i\in \mathcal{V}$
\begin{equation} \label{equ0}
\begin{array}{lll}
& E[\omega_{i,k} | \mathcal{F}_{k-1}'] = \sum_{j=1}^n a_{ij,k} E [ \omega_{ij,k} | \mathcal{F}_{k-1}']  =\mathbf{0}.
\end{array} \end{equation} Similarly, by  Assumption \ref{ass-noise}-b it is shown that 
\begin{align}\label{estimate0} 
& E[\zeta_{i,k} | \mathcal{F}_{k-1}'] = \mathbf{0} ~ ~\forall i\in \mathcal{V}.  \end{align} 
 Then  from \eqref{equ0} \eqref{estimate0}   and  Assumption \ref{ass-noise}-c,   by   \eqref{def-e2} we derive
\begin{align}
& E[ e_{2,k} | \mathcal{F}_{k-1}']=E[ \omega_k | \mathcal{F}_{k-1}']+E[  \zeta_k | \mathcal{F}_{k-1}']+E[v_k | \mathcal{F}_{k-1}']= \mathbf{0}.
\label{equ3}
\end{align}
Since $\mathcal{F}_{k-1} \subset \mathcal{F}_{k-1}^{'}$, by \eqref{equ0} \eqref{equ3}  we see 
\begin{align}
&  E[ \omega_{k} | \mathcal{F}_{k-1}]=  E \big[ E[ \omega_{k} | \mathcal{F}_{k-1}']  \big | \mathcal{F}_{k-1} \big]=  \mathbf{0}~a.s., \label{equ4}\\&E[ e_{2,k} | \mathcal{F}_{k-1}]=  E \big[ E[ e_{2,k} | \mathcal{F}_{k-1}']  \big | \mathcal{F}_{k-1} \big]=  \mathbf{0}~a.s.   \nonumber
\end{align} 
  
Since $a_{ij,k}$ is   adapted to  $\mathcal{F}_{k-1}'$ by \eqref{def-algebras}, from  Assumption \ref{ass-noise}-a it follows that 
 \begin{equation} 
\begin{split}
  &E[ \|   a_{ij,k}   \omega_{ij,k} \|^2   | \mathcal{F}_{k-1}']  \leq     a_{ij,k}^2  E[ \|     \omega_{ij,k} \|^2   | \mathcal{F}_{k-1}']     \leq  a_{ij,k}^2 \mu^2. \nonumber
\end{split} 
\end{equation}
 Since $a_{ij,k}$ is    
independent of  $\mathcal{F}_{k-1}  $ by Assumption \ref{ass-graph}-d,     from  $\mathcal{F}_{k-1} \subset \mathcal{F}_{k-1}^{'}$  by Assumption \ref{ass-graph}-c  we obtain 
 \begin{equation} 
\begin{split}
 E[ \|   a_{ij,k}   \omega_{ij,k} \|^2   | \mathcal{F}_{k-1}] &= E \big[ E[ \|   a_{ij,k}   \omega_{ij,k} \|^2   | \mathcal{F}_{k-1}']  \big| \mathcal{F}_{k-1} \big] 
\\&  \leq     \mu^2    E [ a_{ij,k} ^2 ]  \leq   \mu^2 \eta^2 ~~ \forall i \in \mathcal{V}~a.s.\nonumber
\end{split} 
\end{equation} 
 Then by the   conditional  Minkowski  inequality  $
\big(E[\| \sum_{i=1}^kX_i \|^2   \big | \mathcal{F}] \big)^{\frac{1}{2}} \leq \sum_{i=1}^k \big(E[\| X_i \|^2 \big | \mathcal{F}] \big)^{\frac{1}{2}} $,  and  by    $\omega_{i,k}= \sum_{j=1}^n a_{ij,k}   \omega_{ij,k} $    we derive 
\begin{equation} 
\begin{split} 
\big(E[ \|  \omega_{i,k}  \|^2   \big |  \mathcal{F}_{k-1}  ] \big)^{\frac{1}{2}} & \leq   \sum_{j=1}^n  \big( E[ \|   a_{ij,k}   \omega_{ij,k} \|^2   | \mathcal{F}_{k-1}] \big)^{\frac{1}{2}}   \leq n \mu \eta ~~ \forall i \in \mathcal{V}~a.s. \nonumber
\end{split}
\end{equation}

Similarly, by  Assumption \ref{ass-noise}-b we  derive  
$$ \big(  E[\| \zeta_{i,k} \|^2  | \mathcal{F}_{k-1} ] \big)^{\frac{1}{2}}   \leq n \mu \eta~~ \forall i \in \mathcal{V}~~a.s.$$
Then by the definitions of   $\omega_k  $ and $\zeta_k$  we conclude that 
\begin{equation}\label{n1}
\begin{split}
&E [\| \omega_k  \|^2 | \mathcal{F}_{k-1} ] = \sum_{i=1}^n E [\| \omega_{i,k}  \|^2 | \mathcal{F}_{k-1} ]  \leq  n^3\mu^2 \eta^2~a.s. , 
\\&E [\| \zeta_k  \|^2 | \mathcal{F}_{k-1}]  =\sum_{i=1}^n E [\| \zeta_{i,k}  \|^2 | \mathcal{F}_{k-1} ]\leq  n^3\mu^2 \eta^2~a.s.
   \end{split} 
\end{equation}

By \eqref{gradientnoise} we have
\begin{equation}  
E [\| v_k \|^2 | \mathcal{F}_{k-1}']= \sum_{i=1}^n E [ \| v_{i,k} \|^2 | \mathcal{F}_{k-1}']  \leq c_v(n+ \|X_k\|^2) .\nonumber
\end{equation}
  Then by  noticing that $\mathcal{F}_{k-1} \subset \mathcal{F}_{k-1}^{'}$ and $X_k$ is adapted to $\mathcal{F}_{k-1}$ we have
   \begin{equation} \label{61}
   \begin{split}
E [\| v_k \|^2 | \mathcal{F}_{k-1}]= E\big[ E [\| v_k \|^2 | \mathcal{F}_{k-1}'] \big|\mathcal{F}_{k-1}\big] \leq c_v(n+ \|X_k\|^2) ~~a.s. \end{split}
\end{equation}
Thus,    by \eqref{def-e2} from  \eqref{n1}   \eqref{61}   we obtain 
\begin{equation} 
\begin{split} 
 E[ \| e_{2,k} \|^2 | \mathcal{F}_{k-1} ]  &=  3( E [\| \omega_k  \|^2 | \mathcal{F}_{k-1} ]+E [\| \zeta_k\|^2 | \mathcal{F}_{k-1} ]+E [\| v_k \|^2  | \mathcal{F}_{k-1} ] ) \\ &\leq  6n^3 \mu^2 \eta^2 +3c_v(n+ \|X_k\|^2)  ~a.s. \nonumber 
   \end{split} 
\end{equation}
Hence  \eqref{noise2} holds.

We now  consider properties  of the noise sequence $\{e_{3,k})$  defined  in \eqref{def-e3}.
Since  $X_k $  is  adapted to $\mathcal{F}_{k-1}  $,  by  \eqref{conditional-Laplacian} \eqref{equ4}
we have
\begin{align}
  E[ e_{3,k}| \mathcal{F}_{k-1}] =  E [   \mathcal{\bar{L}}-\mathcal{ L}_k   
 | \mathcal{F}_{k-1}   ]  X_k -E[\omega_k|\mathcal{F}_{k-1}]=\mathbf{0}~a.s. \nonumber 
\end{align}
Since   $X_k ,  \mathcal{L}_k$ are adapted to  $ \mathcal{F}_{k-1}'$ 
  and  $\mathcal{F}_{k-1} \subset\mathcal{F}_{k-1}'$,   by  \eqref{equ0}  we derive 
 \begin{equation}\label{01}
\begin{split}
& E [ \omega_k^T  ( \mathcal{L}_k -\mathcal{\bar{L}})  X_k | \mathcal{F}_{k-1}]  =
  E \Big [   E \big [ \omega_k^T | \mathcal{F}_{k-1}' \big]   ( \mathcal{L}_k -\mathcal{\bar{L}})  X_k  \big | \mathcal{F}_{k-1}\Big] =0~a.s.\end{split}
\end{equation}
 Hence  by \eqref{n1} and Assumption \ref{ass-graph}-d we conclude that
 \begin{equation}\label{01}
\begin{split}
 E [ \| e_{3,k} \|^2 | \mathcal{F}_{k-1}]&=  E[ \|( \mathcal{L}_k -\mathcal{\bar{L}})   X_k \|^2 |\mathcal{F}_{k-1}]   + E [\|   \omega_k \|^2 | \mathcal{F}_{k-1} ]
 + 2 E [ \omega_k^T ( \mathcal{L}_k -\mathcal{\bar{L}})   X_k | \mathcal{F}_{k-1}]  \\& \leq  E [\|  \mathcal{L}_k -\mathcal{\bar{L}}\|^2]  \| X_k \|^2 + n^3\mu^2 \eta^2~a.s.\end{split}
\end{equation} 
Therefore,  \eqref{noise3}  holds.
\hfill  $\blacksquare$

\begin{lem} \label{lemma2}
Let    Assumptions \ref{ass-graph}-a,  \ref{ass-graph}-c,   \ref{ass-graph}-d  and  \ref{ass-noise} hold. Then for any $X \in \Omega$ and $\Lambda \in  \mathbb{R}^{mn}$   
  \begin{equation}  \label{result1} 
  \begin{split}
   E [ \| X_{k+1}- X \| ^2 | \mathcal{F}_{k-1} ] &  \leq \| X_k-X \|^2 +\gamma_k  ^2 \|     \nabla \widetilde{f}(X_{k})   +\mathcal{\bar{L}} \big( \Lambda_k +  X_k\big)  \|^2 \\
& +2 \gamma_k \big(\Phi(X,\Lambda_k) - \Phi(X_k, \Lambda_k) \big) 
 -2\gamma_k \big ( X_k-X \big)^T  \mathcal{\bar{L}}  X_k \\&
+ C_{01}\gamma_k^2 \|  \Lambda(k) +  X_k\|^2 + 3c_v \gamma_k^2\| X_k\|^2+  C_{02}\gamma_k^2 ~~a.s., 
\end{split}
\end{equation}
and
\begin{equation}   \label{result2}
\begin{split}
 E [ \| \Lambda_{k+1}- \Lambda \| ^2  | \mathcal{F}_{k-1}]   &\leq  \| \Lambda_{k} -\Lambda \|^2  +\gamma_k ^2\|  \mathcal{\bar{L}}X_k \|^2   + C_{03}\gamma_k ^2 \\&+2\gamma_k \big( \Phi(X_k,\Lambda_k) -\Phi(X_k, \Lambda )\big) + C_{01} \gamma_k ^2 \| X_k \|^2~~a.s.
\end{split}
\end{equation}
\end{lem}

{\bf Proof:}
By using the non-expansive property  \eqref{pro} of the projection operator,   from \eqref{compact}  we obtain   
 \begin{equation}\label{inequ1}
\begin{split}
  \| X_{k+1}- X \| ^2 & \leq  \big \| X_k -  \gamma_k    \nabla \widetilde{f} (X_k)   -\gamma_k\mathcal{\bar{L}} \big( \Lambda_k +  X_k \big) -X + \gamma_k  \big(e_{1,k}+e_{2,k}\big)  \big \|^2 \\
 & \leq I_0(k) +  \gamma_k ^2 I_1(k)   + 2 \gamma_k I_2(k)  ~~ \forall X \in \Omega,
 \end{split}
\end{equation}
where  $I_0(k)= \| X_k -  \gamma_k   \nabla \widetilde{f} (X_k)   -\gamma_k\mathcal{\bar{L}} \big( \Lambda_k +  X_k \big) -X \|^2$,
$I_1(k)= \| e_{1,k}+e_{2,k}  \|^2$, $I_2(k)=\big(e_{1,k}+e_{2,k}\big)^T
\Big( X_k -  \gamma_k   \nabla \widetilde{f} (X_k)  -\gamma_k\mathcal{\bar{L}} \big( \Lambda_k  +  X_k  \big) -X\Big) $.

    Since  $e_{1,k}$ is adapted to $\mathcal{F}_{k-1}'$ by  \eqref{def-algebras} \eqref{def-e1}, 
   by $\mathcal{F}_{k} \subset \mathcal{F}_{k}'$  and  \eqref{equ3}   we see that 
\begin{equation}\label{equ5}
\begin{split}& E [ e_{1,k}^T e_{2,k}| \mathcal{F}_{k-1}]= E \big[  E [  e_{1,k}^T e_{2,k}   | \mathcal{F}_{k-1}'] \big | \mathcal{F}_{k-1} \big ]  = E \big[  e_{1,k}^T E [  e_{2,k}  | \mathcal{F}_{k-1}'] \big | \mathcal{F}_{k-1} \big ]=0~a.s.
\end{split}
\end{equation} 
 Thus, from  here by \eqref{noise1} and \eqref{noise2}  we derive 
\begin{equation}\label{compon1}
\begin{split}
 E [I_1(k) |\mathcal{F}_{k-1}] &= E [ \| e_{1,k}\| ^2 |\mathcal{F}_{k-1}]+E [  \| e_{2,k}\| ^2 |\mathcal{F}_{k-1}] +
 2E [e_{1,k}  ^T e_{2,k} |\mathcal{F}_{k-1}]
 \\&  \leq C_{01} \|  \Lambda_k +  X_k\|^2 + C_{02}  + 3c_v \| X_k\|^2~a.s.
\end{split}
\end{equation}

Since  $X_k,\Lambda_k$ are adapted to $\mathcal{F}_{k-1}$, by   \eqref{noise1}   \eqref{noise2}  we derive
\begin{equation}\label{compon2}
\begin{split}
&E [I_2(k) |\mathcal{F}_{k-1}] = E \big[  e_{1,k}+e_{1,k}  |\mathcal{F}_{k-1} \big]^T  \big( X_k  -  \gamma_k    \nabla \widetilde{f}(X_k) -\gamma_k\mathcal{\bar{L}}  ( \Lambda_k +  X_k  ) -X\big)=0~a.s.
\end{split}
\end{equation}

Since $I_0(k)$ is adapted to $\mathcal{F}_{k-1}$, combining \eqref{inequ1}, \eqref{compon1}, \eqref{compon2} we   obtain 
\begin{equation} \label{bound01}
\begin{split}
 E [ \| X_{k+1}- X \| ^2 | \mathcal{F}_{k-1}]  &  \leq  I_0(k)  +C_{01}   \gamma_k ^2\|  \Lambda_k +  X_k\|^2 
 +3 c_v \gamma_k^2\| X_k\|^2+ C_{02} \gamma_k ^2 ~~a.s .
    \end{split}
\end{equation}

Note  that
\begin{equation}\label{bound11}
\begin{split}
    I_0(k) & = \| X_k -  \gamma_k  \nabla \widetilde{f}(X_k)  
    -\gamma_k\mathcal{\bar{L}} \big( \Lambda_k+  X_k \big) -X \|^2\\
& \leq  \| X_k-X \|^2+  \gamma_k  ^2 \|    \nabla \widetilde{f}(X_k)  +\mathcal{\bar{L}} \big( \Lambda_k +  X_k \big)  \|^2\\
& - 2 \gamma_k\big ( X_k-X \big)^T\Big(  \nabla \widetilde{f}(X_k) +\mathcal{\bar{L}} \big( \Lambda_k +  X_k \big) \Big).
\end{split}
\end{equation}
Since   $\Phi(X, \Lambda_k) $ is convex in $X\in \Omega$,  by \eqref{gradient}
 we derive 
$$ \Phi(X,\Lambda_k) \geq \Phi(X_k, \Lambda_k)+   (X-X_k)^T \big (   \nabla \widetilde{f}(X_k)+  \mathcal{\bar{L}} \Lambda_k \big),$$
and hence   
$$ -(X_k-X) ^T \big (   \nabla \widetilde{f}(X_k)+  \mathcal{\bar{L}} \Lambda_k\big)  \leq  \Phi(X,\Lambda_k) - \Phi(X_k, \Lambda_k).$$
Then  by  \eqref{bound11} we conclude that 
 \begin{equation} 
 \begin{split}
  I_0(k)  & \leq  \| X_k-X \|^2+  \gamma_k  ^2 \|   \nabla \widetilde{f}(X_k)  +\mathcal{\bar{L}} \big( \Lambda_k +  X_k \big)  \|^2\\
& +2 \gamma_k \big(\Phi(X,\Lambda_k) - \Phi(X_k, \Lambda_k) \big)  
 -2\gamma_k \big ( X_k -X \big)^T  \mathcal{\bar{L}}  X_k , \nonumber
\end{split}
\end{equation}
 which incorporating with   \eqref{bound01}  yields \eqref{result1}.

For any $\Lambda \in \mathbb{R}^n$
 \begin{equation}\label{inequ2}
\begin{split}
   \| \Lambda_{k+1}- \Lambda \| ^2&  =  \| \Lambda_k+\gamma_k \mathcal{ \bar{L}}  X_k-\Lambda +\gamma_ke_{3,k} \|^2 \\ & = I_3(k)  +\gamma_k ^2\| e_{3,k} \|^2  +2\gamma_k  e_{3,k}^T\big(\Lambda_{k}+\gamma_k \mathcal{ \bar{L}} X_k -\Lambda\big),\end{split}
\end{equation}
where $I_3(k)=\| \Lambda_k+\gamma_k \mathcal{ \bar{L}}   X_k -\Lambda \|^2 .$

Since  $X_k$  and $\Lambda_k$ are  adapted to $\mathcal{F}_{k-1}$,
  from \eqref{noise3}   we see 
  \begin{equation} 
\begin{split}
 & E[ e_{3,k}^T(\Lambda_{k}+\gamma_k \mathcal{ \bar{L}}  X_k -\Lambda) | \mathcal{F}_{k-1}]  =
E[ e_{3,k}^T | \mathcal{F}_{k-1}] (\Lambda_{k}+\gamma_k \mathcal{ \bar{L}}  X_k -\Lambda)=0~a.s. \nonumber
\end{split}
\end{equation}
 Noticing that  $I_3(k)$ is adapted to $ \mathcal{F}_{k-1}$,    from here  by  \eqref{noise3} \eqref{inequ2} we obtain 
 \begin{equation}\label{bound02}
\begin{split}
 E[   \| \Lambda_{k+1}- \Lambda \| ^2  | \mathcal{F}_{k-1}] \leq  I_3(k)
  + C_{01} \gamma_k ^2 \| X_k \|^2 + C_{03}\gamma_k ^2~~a.s .
  \end{split}
\end{equation}

By the definition of    $\Phi(X, \Lambda) $,  we derive 
$$ \Phi(X_k, \Lambda_k) =\Phi(X_k, \Lambda )+   (\Lambda_k-\Lambda )^T  \mathcal{\bar{L}}  X_k$$
and hence
\begin{equation} 
\begin{split}
I_3(k) &=
\| \Lambda_{k} -\Lambda \|^2  +\|  \gamma_k \mathcal{\bar{L}}X_k \|^2+2\gamma_k(\Lambda_k -\Lambda )^T \mathcal{\bar{L}}X_k \\& =\| \Lambda_{k} -\Lambda \|^2  +\|  \gamma_k \mathcal{\bar{L}}X_k  \|^2  +2\gamma_k \big( \Phi(X_k ,\Lambda_k) -\Phi(X_k , \Lambda )\big), \nonumber
\end{split}
\end{equation}
which incorporating  with \eqref{bound02} yields   \eqref{result2}.
\hfill  $\blacksquare$
 
\textbf{Proof of Theorem  \ref{thm2}:}  
Summing up both sides of \eqref{result1} and  \eqref{result2},   and by replacing  $(X, \Lambda)$  with 
  $(X^{*}, \Lambda^{*})$  we obtain 
\begin{equation} \label{sum1}
\begin{split}
 & E [ \| X_{k+1}- X^{*} \| ^2 | \mathcal{F}_{k-1} ]  +E [ \| \Lambda_{k+1}- \Lambda^{*} \| ^2  | \mathcal{F}_{k-1}]   
 \\  &\leq  \| X_k -X^{*} \|^2+  \| \Lambda_{k} -\Lambda^{*} \|^2   +  \gamma_k  ^2 \|   \nabla \widetilde{f}(X_k)  +\mathcal{\bar{L}} \big( \Lambda_k +  X_k \big)  \|^2
  \\& +2 \gamma_k \big(\Phi(X^{*},\Lambda_k) - \Phi(X_k, \Lambda^{*}) \big)  
-2\gamma_k \big ( X_k-X^{*} \big)^T  \mathcal{\bar{L}}  X_k    \\&
+ C_{01} \gamma_k^2\|  \Lambda_k +  X_k\|^2 +\gamma_k^2\|   \mathcal{\bar{L}}X_k \|^2   + 
  \gamma_k ^2 ( C_{01}+3c_v) \| X_k \|^2 + ( C_{02}+ C_{03}) \gamma_k ^2 ~a.s.
   \end{split}
\end{equation}

Since $(X^{*}, \Lambda^{*})$ is   a saddle point for $\Phi(X,\Lambda)$,  by Lemma \ref{lem3} $X^{*}$ is the optimal solution to  the problem \eqref{problem2}. Then from  Lemma \ref{lem2}  it follows that 
\begin{equation}\label{consensus}
\mathcal{\bar{L}} X^{*}=\mathbf{0},~\textrm{and }~\mathcal{\bar{L}}   X_k= \mathcal{\bar{L}} (   X_k -X^{*}),
\end{equation} 
and hence 
  \begin{equation} 
\begin{split}
&   \nabla \widetilde{f}(X_k) + \mathcal{\bar{L}} (  \Lambda_k+X_k )  =   \nabla \widetilde{f}(X_k) -  \nabla \widetilde{f}(X^{*})     +
  \mathcal{\bar{L}} (  \Lambda_k-\Lambda^{*} )+ \mathcal{\bar{L}}(  X_k-X^{*})
 + \mathcal{\bar{L}}\Lambda^{*} +  \nabla \widetilde{f}(X^{*}).\nonumber
  \end{split}
\end{equation} 
Then by \eqref{Lipschitz}  we obtain 
 \begin{equation} \label{b1}
\begin{split}
  &\| \nabla \widetilde{f}(X_k) + \mathcal{\bar{L}} (  \Lambda_k+X_k )  \|^2    \\&\leq 
4   \big(   \|  \nabla \widetilde{f}(X_k)  - \nabla \widetilde{f} (X^{*})\| ^2+  \| \mathcal{\bar{L}} (  \Lambda_k -\Lambda^{*} ) || ^2  +  \| \mathcal{\bar{L}} ( X_k-X^{*} ) || ^2 +  \|\mathcal{\bar{L}}\Lambda^{*} +\nabla \widetilde{f}(X^{*})  \|^2  \big)
  \\& \leq  4c_1 \|  \Lambda_k  -\Lambda^{*}\| ^2+ (4c_1+ 4L_f^2) \|X_k-X^{*}\|^2+ c_2,  
  \end{split}
\end{equation} 
where $c_1=  \|  \mathcal{\bar{L} } \|^2$, $c_2= 4\|\mathcal{\bar{L}}\Lambda^{*} +\nabla \widetilde{f}(X^{*})  \|^2 $.
From \eqref{consensus} we derive
\begin{align}
& \|\mathcal{\bar{L}}   X(k) \|^2  \leq c_1 \| X_k- X^{*} \|^2.\label{b2} \end{align}  
Note that $ \|  \Lambda_k+  X_k\|^2   \leq 3(  \|  \Lambda_k -\Lambda^{*}\|^2 +\|   X_k -X^{*}\|^2+\|  \Lambda^{*} +  X^{*}\|^2)$ and  $\| X_k \|^2 \leq 2(  \| X_k-X^{*} \|^2+ \|  X^{*} \|^2 ).$
Then by   \eqref{sum1}, \eqref{b1} and  \eqref{b2} we derive
\begin{equation}   \label{sum2}
\begin{split}
 & E [ \| X_{k+1}- X^{*} \| ^2 | \mathcal{F}_{k-1} ]  +E [ \| \Lambda_{k+1}- \Lambda^{*} \| ^2  | \mathcal{F}_{k-1}]   
 \\  &\leq  \big( 1+  (5 c_1+ 5 C_{01}+ 4 L_f^2+6c_v) \gamma_k^2 \big)\| X_k -X^{*} \|^2 \\&+   \big( 1+(4c_1+ 3C_{01} ) \gamma_k^2 \big)\| \Lambda_{k} -\Lambda^{*} \|^2   +\gamma_k ^2  \big( c_2+C_{02}+ C_{03} +3C_{01} \| \Lambda^{*} +X^{*}\|^2\\&+ 2(3 c_v+ C_{01} )\| X^{* }\|^2) 
 +2 \gamma_k \big(\Phi(X^{*},\Lambda_k ) - \Phi(X_k, \Lambda^{*}) \big)  -2 \gamma_k \big ( X_k-X^{*} \big)^T  \mathcal{\bar{L}}  X_k ~a.s.  
   \end{split}
\end{equation}

Since $\mathcal{\bar{L} } $ 
is the  Laplacian matrix of  some  connected undirected graph by Assumption   \ref{ass-graph}-b, 
from \eqref{consensus}  and Lemma \ref{lem1} it follows that \begin{equation}\label{b5}
\begin{split}
& \big ( X_k-X^{*} \big)^T  \mathcal{\bar{L}}  X_k   =\big ( X_k-X^{*} \big)^T  \mathcal{\bar{L}} \big(  X_k - X^{*}\big)    \geq 0.
\end{split}
\end{equation}  
Noticing $X_k \in \Omega$,  by  definition of  the  saddle point we see 
 \begin{equation} 
\Phi(X^{*} , \Lambda_k) \leq  \Phi(X^{*},   \Lambda^{*})  \leq \Phi(X_k , \Lambda^{*}) ~~ \forall k \geq 0. \nonumber
\end{equation}
Then by setting $V_{k}= \| X_k-X^{*} \|^2 +   \| \Lambda_{k} -\Lambda^{*} \|^2$, from   \eqref{sum2}  \eqref{b5}  we derive 
\begin{equation} \begin{split}
 & E\big  [ V_{k+1} | \mathcal{F}_{k-1} \big ]   \leq  ( 1 + C_{11}\gamma_k^2  ) V_k+ C_{12}\gamma_k ^2~~a.s. ,\nonumber    \end{split}
\end{equation}
where $C_{11}=5 c_1+ 5 C_{01}+ 4 L_f^2+ 6c_v, $ and   $C_{12}=c_2+C_{02}+ C_{03} +3C_{01} \| \Lambda^{*} +X^{*}\|^2+ 2(3c_v+ C_{01} ) \| X^{* }\|^2 $.

Consequently,  by  Assumption \ref{ass-stepsize}  and  Lemma \ref{martingale}  in Appendix we conclude that  $\| X_{k}- X^{*} \| ^2 +   \| \Lambda_{k}- \Lambda^{*} \| ^2$ converges a.s.   \hfill   $\blacksquare$
 
\subsection{Consensus and Consistency}
The following theorem   shows  that the estimates    given by all agents  reach a consensus belonging  to the optimal solution set  of problem \eqref{problem1}.  
 \begin{thm}\label{thm1}    Let $\{ x_{i,k}\}$ and  $\{ \lambda_{i,k}\}$   be produced by the algorithm  \eqref{algorithm1}
with any   initial values $x_{i,0},~ \lambda_{i,0}$.   
 Let Assumptions \ref{ass-set},  \ref{ass-graph},    \ref{ass-noise}, and      \ref{ass-stepsize}  hold.  Then 

i) (Consensus) $\lim\limits_{k \rightarrow \infty} (x_{i,k} -x_{j,k})=0~~\forall i ,j \in \mathcal{V}~~a.s.$

ii) (Consistency) 
 \begin{equation}\label{mainresults}
  \lim_{k \rightarrow \infty}  d(x_{i,k} , \Omega_o^{*})=0  ~~\forall i \in \mathcal{V} ~~a.s., 
\end{equation}
where  $d(X,A)=\inf\limits_{\theta \in A} \| \theta-X\|. $
Moreover, if $f(\cdot)$ has a   unique optimal solution  $x^{*}$,  then
\begin{equation}\label{mainresults2}
  \lim_{k \rightarrow \infty}  x_{i,k} =x^* ~~\forall i \in \mathcal{V} ~~a.s. 
  \end{equation}
  \end{thm}

\textbf{Proof:}  By setting
\begin{equation}\label{ }
\begin{split}
& \theta=\begin{pmatrix}
      X   \\
      \Lambda
\end{pmatrix},~~e_k=\begin{pmatrix}
      e_{1,k}+e_{2,k}    \\
      e_{3,k}  
\end{pmatrix} , ~~\Phi=  \Omega \times \mathbb{R}^{mn},  \\
& g(\theta)= \begin{pmatrix}
     - \nabla \widetilde{f}(X) -\mathcal{\bar{L}}\big(X+\Lambda  \big)   \\
      \mathcal{\bar{L}} X
\end{pmatrix} , \nonumber
\end{split}
\end{equation}
 we can rewrite \eqref{compact} in the form of algorithm \eqref{constrained}  with  $Y_k=g(\theta_k)+e_k$.
We intend to use   Lemma \ref{CSA}   in Appendix to prove the theorem.  Thus, we have to verify   B1-B4.  

 Since $X_k,~\Lambda_k$ are bounded  a.s., from \eqref{noise1}, \eqref{noise2}, \eqref{noise3} we conclude that 
\begin{equation}
\begin{split}
  &E[ \| e_{k}\|^2 ]   \leq 2E[ \| e_{1,k}\|^2 ] +2E[ \| e_{2,k}\|^2  ] +E[ \| e_{3,k}\|^2 ] < \infty . \nonumber
\end{split}
\end{equation}
Thus, $E [ \|Y_k\|^2] =2E [  \|g(\theta_k)||^2] +2E[ \| e_{k}\|^2]  < \infty , $ and hence B1 holds.
From     \eqref{noise1}  \eqref{noise2}  \eqref{noise3}   it follows that   B2 holds.
By the definition of $g(\theta)$, from Assumption \ref{ass-set}-b it is seen  that   B3 holds.
By Theorem \ref{thm2} we conclude that  $\theta_k$ is bounded a.s., and hence B4 holds.

 In summary, we have validated B1-B4. Then  by    Lemma \ref{CSA}  we conclude that  
 $(X_k, \Lambda_k) $ converge a.s.   to  some limit   set of the following  projected ODE in $\Omega \times \mathbb{R}^{mn}$: 
 \begin{equation} \label{UDI}
\begin{split}
& \dot{X}(t)=  -\nabla \widetilde{f}(X(t)) -\mathcal{\bar{L}}\big(X(t)+\Lambda(t)  \big)-Z(t) , ~Z(t) \in  N_{\Omega}(X(t)),\\&
\dot{ \Lambda}(t)=\mathcal{\bar{L}} X(t),
\end{split}
\end{equation}
where  $Z(\cdot)  $ is the projection or constraint  term,   the minimum force needed  to keep $X(\cdot)$ in $\Phi.$

  Define $V(X,\Lambda) =\| X-X^{*}\|^2+ \| \Lambda-\Lambda^{*}\|^2$. By \eqref{UDI} we derive
  \begin{equation}
\begin{array}{ll}
\dot{V}(X,\Lambda)   
&=( X-X^{*})^T \dot{X}+(\Lambda-\Lambda^{*})^T \dot{\Lambda}
 \\&=-( X-X^{*})^T \widetilde{f}(X) -( X-X^{*})^T\mathcal{\bar{L}}\big(X +\Lambda   \big)  -( X-X^{*})^T Z 
+(  \Lambda-\Lambda^{*})^T \mathcal{\bar{L}} X. \nonumber \end{array}
\end{equation}
Since   $Z(t) \in  N_{\Omega}(X(t))$,  by the  definition of normal cone we derive
$( X-X^{*})^T Z  \geq 0.$   Since  $ \mathcal{\bar{L}} $ is symmetric, by \eqref{consensus} we  derive 
  \begin{equation} \label{ineqs}
\begin{array}{ll}
\dot{V}(X,\Lambda)  & \leq -( X-X^{*})^T \widetilde{f}(X)-  X ^T\mathcal{\bar{L}} X -X ^T\mathcal{\bar{L}}    \Lambda   +  
  \Lambda^T  \mathcal{\bar{L}} X  -(\Lambda^{*})^T \mathcal{\bar{L}} (X-X^{*})  \\&  \leq -( X-X^{*})^T \big( \widetilde{f}(X)+\mathcal{\bar{L}} \Lambda^{*}\big)- X^T\mathcal{\bar{L}} X  \\& \leq  \Phi(X^{*}, \Lambda^{*}) -\Phi(X, \Lambda^{*}) -X^T\mathcal{\bar{L}} X,  
   \end{array}
\end{equation}
where in the last inequality  we have used   $\Phi(X^{*},\Lambda^{*}) \geq \Phi(X, \Lambda^{*})+   (X^{*}-X )^T \big (   \nabla \widetilde{f}(X )+  \mathcal{\bar{L}} \Lambda^{*} \big)$ since $\Phi(X, \Lambda^{*})$ is convex with respect to $X.$ Noting that   $\mathcal{\bar{L}} $ is positive semi-definite, by the definition of saddle point we derive
  \begin{equation}\label{dot}
\begin{array}{ll}
\dot{V}(X,\Lambda)  & \leq0. \nonumber  \end{array}
\end{equation}

By the LaSalle invariant  theorem \cite{nonlinear},  the trajectories produced by \eqref{UDI} converge to the largest  invariant set contained  in the set   $S=\{(X,\Lambda)\in \Omega \times \mathbb{R}^{mn}: \dot{V}(X,\Lambda) =0\}.$
By \eqref{ineqs} it is clear that  $S=\{(X,\Lambda)\in \Omega \times \mathbb{R}^{mn}:X^T\mathcal{\bar{L}} X=0,
~\Phi(X^{*}, \Lambda^{*}) -\Phi(X, \Lambda^{*}) =0\}.$ 

If $X^T\mathcal{\bar{L}} X=0$,   then by noticing that  $\mathcal{\bar{L}}$ is the Laplacian matrix of an undirected  connected  graph, from  Lemma \ref{lem1}  we have  $X=  \mathbf{1}  \otimes  x$ for some $x \in \mathbb{R}^m$.
Since $(X^{*}, \Lambda^{*})$ is a saddle point of $\Phi(X,\Lambda)$,   $X^{*}$ is an optimal solution to the problem \eqref{problem2} by Lemma \ref{lem3}. Thus,   $\Phi(X^{*}, \Lambda^{*})=\widetilde{f}(X^{*})+(\Lambda^{*})^T \mathcal{\bar{L}} X^{*}=\widetilde{f}(X^{*})=f^{*}.$  If $ \Phi(X, \Lambda^{*})-\Phi(X^{*}, \Lambda^{*}) =0$, then  from  $X=  \mathbf{1}  \otimes  x$ and $X\in \Omega$ we conclude that  $f(x)=f^*,~ x\in \Omega_{o}.$
Thus, $x$ is also an  optimal  solution to problem \eqref{problem1}, and hence 
$S=\{(X,\Lambda):X=  \mathbf{1}  \otimes  x, ~ x \in \Omega_o^{*}, ~\Lambda \in \mathbb{R}^{mn}\}.$ Therefore,  $(X_k, \Lambda_k) $ converges to   the largest  invariant set in set $S$.  Consequently,  the estimates given by all agents finally reach consensus,  and  hence  \eqref{mainresults}  holds. 

Furthermore,   if $\Omega_o^{*}=\{x^*\}$, then by   \eqref{mainresults} we derive  \eqref{mainresults2}. 
The proof is completed.  
  \hfill $\blacksquare$

  \section{Asymptotic Properties}\label{sec:Normality}
  In this section, we establish  asymptotic properties of the distributed  primal-dual algorithm \eqref{algorithm1} when there is no constraint, i.e.,  $\Omega_i =\mathbb{R}^m~~\forall i \in \mathcal{V}$. 
\subsection{Dimensionality Reduction }
We  now  introduce a linear transformation to  the algorithm \eqref{compact}.
Note that    $\mathcal{\bar{L}}$  is   the Laplacian matrix of an undirected  connected  graph by Assumption \ref{ass-graph}-b. Then  by Lemma  \ref{lem1}  
$\mathcal{\bar{L}}$  has a simple zero eigenvalue while all other eigenvalues are positive. 
Thus,  there exists an  orthogonal matrix $\mathcal{V}=(\mathcal{V}_1~\mathcal{V}_2)$, where $\mathcal{V}_2=\frac{\mathbf{1}}{\sqrt{n}}$ and each column of $ \mathcal{V}_1 \in \mathbb{R}^{n \times (n-1)}$ is an eigenvector corresponding  to   some  positive eigenvalue of  $\mathcal{\bar{L}}$,  such that \begin{equation}\label{diagnoal}
 \mathcal{V}^T   \mathcal{ \bar{L}}  \mathcal{V}=\begin{pmatrix}
  \mathcal{ S}   & \mathbf{ 0 }  \\
     \mathbf{0}  &  0
\end{pmatrix}
\end{equation}
where      $\mathcal{S} =diag\{\kappa_2,  \cdots, \kappa_n\} \in \mathbb{R}^{(n-1) \times (n-1)}$ with $\kappa_i,i=2,\cdots,n$   being  positive eigenvalues of  $\mathcal{\bar{L}}$.

By multiplying both sides of \eqref{diagnoal} from left with $\mathcal{V}$, 
it follows that \begin{align}\label{left}
 \mathcal{ \bar{L}}   \mathcal{V}=(\mathcal{V}_1~\mathcal{V}_2)\begin{pmatrix}
  \mathcal{ S}   & \mathbf{ 0 }  \\
     \mathbf{0}  &  0
\end{pmatrix} =(\mathcal{V}_1\mathcal{S}~\mathbf{0} ).
\end{align}
Similarly, by multiplying both sides of \eqref{diagnoal} from  right with $\mathcal{V}^T$,  
we obtain \begin{align}\label{right}
 \mathcal{V}^T\mathcal{\bar{L}}=\begin{pmatrix}
  \mathcal{ S}   & \mathbf{ 0 }  \\
     \mathbf{0}   &  0
\end{pmatrix}  \begin{pmatrix}
      \mathcal{V}_1^T   \\
   \mathcal{   V}_2^T  
\end{pmatrix}= \begin{pmatrix}
     \mathcal{S} \mathcal{V}_1^T   \\
   \mathbf{0} 
\end{pmatrix}.
\end{align}

 Let $(X^{*}, \Lambda^{*})  $  be    the primal-dual  solution pair of the  problems  \eqref{problem2} and \eqref{pdual} when $ \Omega_i= \mathbb{R}^{m}~\forall i \in \mathcal{V}$. Then by Lemma  \ref{lem3},   $(X^{*}, \Lambda^{*})  $ satisfies    
 $$\Phi(X^{*} , \Lambda) \leq  \Phi(X^{*},   \Lambda^{*})  \leq \Phi(X , \Lambda^{*})~~\forall X, \Lambda \in \mathbb{R}^{mn},$$
 and hence  \begin{equation}\label{optimalcd}
 \nabla \widetilde{f}(X^*)+ ( \mathcal{ \bar{L}}  \otimes \mathbf{I}_m) \Lambda^*=\mathbf{0} ,~~( \mathcal{ \bar{L}}  \otimes \mathbf{I}_m)X^*=\mathbf{0} .
\end{equation}
The first equality in \eqref{optimalcd} is   the optimality condition for $\min\limits_X \widetilde{f}(X)+X^T ( \mathcal{ \bar{L}}  \otimes \mathbf{I}_m) \Lambda^* $, where   the   minimum is attained at $X^*$.

Therefore,  from \eqref{compact} \eqref{optimalcd} and $\Omega=\mathbb{R}^{mn}$ it follows that 
 \begin{align}  X_{k+1}-X^{*}&=   X_k -X^{*}-  \gamma_k   \big(  \nabla \widetilde{f}(X_{k}) - \nabla \widetilde{f}(X^*) \big)  \nonumber \\& -\gamma_k  ( \mathcal{ \bar{L}}  \otimes \mathbf{I}_m)   \big( \Lambda_k-\Lambda^* +  X_k -X^*\big)  + \gamma_k  \big(e_{1,k}+e_{2,k}\big) ,  \label{primal1} \\
\Lambda_{k+1}-\Lambda^*&=\Lambda_k-\Lambda^* +\gamma_k     ( \mathcal{ \bar{L}}  \otimes \mathbf{I}_m)   (X_k-X^*)+  \gamma_k  e_{3,k}  .    \label{dual1}
\end{align}

Define $$\widetilde{\Lambda}_{1,k}\triangleq (\mathcal{V}_1^T \otimes \mathbf{I}_m) (\Lambda_k-\Lambda^*),~\widetilde{\Lambda}_{2,k}\triangleq(\mathcal{V}_2^T \otimes \mathbf{I}_m)(\Lambda_k-\Lambda^*).$$
Then by multiplying both sides of \eqref{dual1} with $\mathcal{V}^T \otimes \mathbf{I}_m$ from left,     by  the rule of Kronecker product 
\begin{equation}\label{pkrp}
(A \otimes B)(C \otimes D)=(A \otimes C)(B \otimes D)
\end{equation}
 and  by \eqref{right} we obtain  
 \begin{equation}
\begin{split}
\begin{pmatrix}
      \widetilde{\Lambda}_{1,k+1}     \\
      \widetilde{\Lambda}_{2,k+1} 
\end{pmatrix}
= \begin{pmatrix}
      \widetilde{\Lambda}_{1,k}     \\
      \widetilde{\Lambda}_{2,k} 
\end{pmatrix} +\gamma_k    \begin{pmatrix}
     \mathcal{S} \mathcal{V}_1^T \otimes \mathbf{I}_m  \\
   \mathbf{0} 
\end{pmatrix} \widetilde{X}_k +  \gamma_k  \begin{pmatrix}
    \mathcal{V}_1^T  \otimes \mathbf{I}_m   \\
      \mathcal{V}_2^T  \otimes \mathbf{I}_m
\end{pmatrix} e_{3,k}  .   \nonumber
\end{split}
\end{equation}
Hence 
\begin{equation}\label{dual2}
  \widetilde{\Lambda}_{1,k+1} =  \widetilde{\Lambda}_{1,k} +\gamma_k (\mathcal{S} \mathcal{V}_1^T \otimes \mathbf{I}_m) \widetilde{X}_k+ \gamma_k   (\mathcal{V}_1^T \otimes \mathbf{I}_m) e_{3,k}.
\end{equation}

Since  $\mathcal{V}\mathcal{V}^T=\mathbf{I}_n$,  by \eqref{left}  and  by \eqref{pkrp} we derive
\begin{equation} 
\begin{split}
  ( \mathcal{ \bar{L}}  \otimes \mathbf{I}_m)   \big( \Lambda_k-\Lambda^*) &= ( \mathcal{ \bar{L}}  \otimes \mathbf{I}_m) (\mathcal{V}\otimes \mathbf{I}_m) ( \mathcal{V}^T \otimes \mathbf{I}_m) \big( \Lambda_k-\Lambda^*) \\&=(\mathcal{V}_1\mathcal{S}  \otimes \mathbf{I}_m ~\mathbf{0} )\begin{pmatrix}
      \widetilde{\Lambda}_{1,k}    \\
      \widetilde{\Lambda}_{2,k}  
\end{pmatrix} =( \mathcal{V}_1\mathcal{S}  \otimes \mathbf{I}_m) \widetilde{\Lambda}_{1,k} . \nonumber
\end{split}
\end{equation}
Then by setting $\widetilde{X}_k \triangleq X_k-X^*,$    from \eqref{primal1} we derive
\begin{equation} \label{primal2}
\begin{split}\widetilde{ X}_{k+1} =  \widetilde{ X}_k& -  \gamma_k   \big(  \nabla \widetilde{f}(\widetilde{X}_{k}+X^*) - \nabla \widetilde{f}(X^*) \big)  \\& -\gamma_k   ( \mathcal{ \bar{L}}  \otimes \mathbf{I}_m)  \widetilde{X}_k-\gamma_k  (\mathcal{V}_1\mathcal{S}  \otimes \mathbf{I}_m) \widetilde{\Lambda}_{1,k}     + \gamma_k  \big(e_{1,k}+e_{2,k}\big).
 \end{split}
 \end{equation}

\subsection{Asymptotic Normality and Efficiency}  
To investigate  the  asymptotic properties of the  algorithm  \eqref{dual2}\eqref{primal2},   we need the following conditions. 

\begin{ass}\label{ass-stepsize1}  $\gamma_k = \frac{1}{k^\nu}$ with $\nu \in (\frac{2}{3},1).$
\end{ass}
 
\begin{ass}\label{ass-function}

  a)  $f(\cdot)$ is strictly convex and the unique  optimal solution is $x^*$.
 
  b) The Hessian matrix of $f_i(\cdot) $ at point  $x^*$  is $H_i$, and   $   \sum_{i=1}^n H_i$ is  positive definite.
  
  c) There exists   a constant $c>0$ such that   $\| \nabla f_i(x)-\nabla f_i(x^*)-H_i(x-x^*)\| \leq  c \|x-x^*\|^2~\forall i\in \mathcal{V}.$ 
\end{ass}

\begin{rem}
 By Assumption \ref{ass-function}-b,   the   Hessian matrix $\nabla^2f(x^*) $ is  positive definite.
If  in addition,  for any $i\in \mathcal{V}$,  the Hessian matrix function $ \nabla^2 f_i(\cdot) $ is globally Lipschitz, then by \cite[Lemma 1.2.4]{Nesterov} we derive
$$\| \nabla f_i(y)-\nabla f_i(x)- \nabla^2 f_i(x)(y-x)\| \leq   \frac{M}{2} \|y-x\|^2,$$
where   $M>0$ is a constant. Hence Assumption \ref{ass-function}-c holds. 
\end{rem}

\begin{ass} \label{ass-independency}

 a) For any $i \neq j \in \mathcal{V}$, $v_{i,k} $ and $v_{j,k}$ are conditionally  independent  given  $\mathcal{F}_{k-1}'.$

b) For any $(i_1,j_1)\neq (i_2,j_2)$ with $i_1,i_2,j_1, j_2 \in \mathcal{V}$,  
$\omega_{i_1j_1,k}$ and $ \omega_{i_2j_2,k}$ are  conditionally  independent  given  $\mathcal{F}_{k-1}'$,
 $\zeta_{i_1j_1,k}$ and $\zeta_{i_2j_2,k}$  are  conditionally  independent  given  $\mathcal{F}_{k-1}'.$ 

 c)  For any $i , j \in \mathcal{V}$, $v_{i,k}$, $\omega_{ij,k}$, and $\zeta_{ij,k}$ are  conditionally  independent  given   $\mathcal{F}_{k-1}'.$
 
 d) For any $i\in \mathcal{V}$,  $v_{i,k}$ and $\mathcal{L}_k$ are conditionally independent  given $\mathcal{F}_{k-1}.$
\end{ass}

Define \begin{equation}\label{symbol}
\begin{split}  & R_{\omega,i} \triangleq \sum_{j=1}^n \sigma_{ij}R_{\omega,ij},  \quad   R_{\zeta,i} \triangleq \sum_{j=1}^n \sigma_{ij}R_{\zeta,ij} ,~~R_v \triangleq diag \big \{R_{v,1},\cdots,R_{v,n} \big) ,\\ &
R_{\omega}    \triangleq    diag \big \{R_{\omega,1} ,\cdots,R_{\omega,n}  \big\},  ~~R_{\zeta}   \triangleq     diag  \big \{R_{\zeta,1} ,\cdots,R_{\zeta,n}  \big \},  \\&  S_1\triangleq E[\big((  \mathcal{ L} _k -  \mathcal{ \bar{L}}   )\mathcal{V}_1S^{-1} \mathcal{V}_1^T \otimes \mathbf{I}_m \big)\nabla \widetilde{f}(X^*) \nabla \widetilde{f}(X^*)^T\big(\mathcal{V}_1S^{-1} \mathcal{V}_1^T(  \mathcal{ L} _k -  \mathcal{ \bar{L}} )^T \otimes \mathbf{I}_m \big)],\\&\mathcal{H}=diag\{H_1,\cdots, H_n\},~~ S_2=R_{v}+R_{\omega}+R_{\zeta}.
 \end{split}
\end{equation}  
\begin{thm}\label{thm3}  (Asymptotic Normality) Set $\Omega_i=\mathbf{R}^m~\forall i \in \mathcal{V}.$ 
 Let Assumptions \ref{ass-set}-b,   \ref{ass-graph},   \ref{ass-noise2}, \ref{ass-stepsize1},  \ref{ass-function},     and   \ref{ass-independency} hold.    Then  
 $  \theta_k=col\{ \widetilde{X}_k   ,     \widetilde{ \Lambda}_{1,k}\}$  is asymptotically normal: 
$$ \theta_k /\sqrt{\gamma_k} \xlongrightarrow [k \rightarrow \infty]{d}  N(\mathbf{0},\Sigma), $$
 where  $\Sigma=\int_{0}^{\infty} e^{Ft}\Sigma_1 e^{F^Tt} dt,$
\begin{equation}\label{hessian}
F\triangleq - \begin{pmatrix}
      ( \mathcal{ \bar{L}}  \otimes \mathbf{I}_m)  +\mathcal{H} & \mathcal{V}_1\mathcal{S} \otimes \mathbf{I}_m \\
    -\mathcal{S} \mathcal{V}_1^T \otimes \mathbf{I}_m  & 0  
\end{pmatrix} ,
\end{equation}   and
$$ \Sigma_1=\begin{pmatrix}
      & S_1+S_2  & -R_{\omega} (\mathcal{V}_1 \otimes \mathbf{I}_m )\\
      &   - (\mathcal{V}_1^T\otimes \mathbf{I}_m) R_{\omega} &(\mathcal{V}_1^T \otimes \mathbf{I}_m ) R_{\omega} ( \mathcal{V}_1 \otimes \mathbf{I}_m )
\end{pmatrix}.$$ 
  \end{thm}

\begin{thm}\label{thm4} Set $\Omega_i=\mathbf{R}^m~\forall i \in \mathcal{V}.$ 
 Let Assumptions \ref{ass-set}-b,   \ref{ass-graph},  \ref{ass-noise2},   \ref{ass-stepsize1},  \ref{ass-function},     and   \ref{ass-independency} hold. 
 Define $  \bar{\theta}_n=\frac{1}{n}\sum_{k=1}^n\theta_k.$
 Then  $\{\bar{\theta}_k\}$ is  asymptotically  efficient:
$$ \sqrt{k} \bar{\theta}_k \xlongrightarrow [k \rightarrow \infty]{d}  N\big(\mathbf{0}, F^{-1}\Sigma_1(F^{-1})^T \big). $$
\end{thm}

 \subsection{Proof of Theorems \ref{thm3} and \ref{thm4}}
   
    Before proving the results, we give some lemmas to be used in  the proof of Theorem \ref{thm3}.
 \begin{lem} \label{Hurwitz}\cite[Lemma 2]{Sayed_TSP_2015}  Let a block matrix $F$ have the following form
 $$F=-\begin{pmatrix}
      X   &Y^T \\
    -Y  &  \mathbf{0} 
\end{pmatrix},$$
 and let $X \in \mathbb{R}^{p \times p}$ be positive definite and $Y \in \mathbb{R}^{p \times q}$ be of  full row rank.
  Then the matrix $F$  is Hurwitz.
\end{lem}

\begin{lem} \label{lem-hessian}  Let Assumption  \ref{ass-graph}-b and Assumption  \ref{ass-function}-b hold. Then  $F$  defined by \eqref{hessian}  is Hurwitz. 
\end{lem}
 {\bf Proof:}
 Since $H_i~\forall i\in \mathcal{V}$ are Hessian matrices of convex functions, $\mathcal{H} $ is semi-positive definite.
The matrix  $  \mathcal{ \bar{L}}   $ is   semi-positive definite since it is a Laplacian matrix of an undirected graph.  
 Therefore, a nonzero vector $x \in \mathbb{R}^{mn}$ satisfies 
 $x^T(\mathcal{ \bar{L}}  \otimes \mathbf{I}_m+\mathcal{H}) x=0$   if and  only if 
 \begin{equation}\label{iif}
x^T ( \mathcal{ \bar{L}}  \otimes \mathbf{I}_m) x=0,~x^T\mathcal{H} x=0 .
\end{equation}  

Since $  \mathcal{ \bar{L}}   $ is the Laplacian matrix of an undirected connected graph,  by Lemma \ref{lem1} a nonzero vector $x \in \mathbb{R}^{mn}$ satisfies $x^T ( \mathcal{ \bar{L}}  \otimes \mathbf{I}_m) x=0$ if and  only if $x=\mathbf{1} \otimes u~\forall u \neq \mathbf{0} \in \mathbb{R}^m . $ 
However,  by Assumption \ref{ass-function}-b $$(\mathbf{1} \otimes u)^T \mathcal{H} (\mathbf{1} \otimes u)= u^T (\sum_{i=1}^n H_i )u>0~~\forall u \neq \mathbf{0} . $$
Therefore,  the two equalities in \eqref{iif} do not hold  simultaneously.   
Thus,  $ ( \mathcal{ \bar{L}}  \otimes \mathbf{I}_m)  +\mathcal{H}$ is positive definite.

 Note that  $\mathcal{S}\mathcal{V}_1^T$ is of full row rank.  Then,   by Lemma  \ref{Hurwitz} we see that 
 $F$ defined by \eqref{hessian} is Hurwitz. 
 \hfill $\blacksquare$
 
\textbf{ Proof of Theorem \ref{thm3}:}  
Set
\begin{equation}\label{gtheta}
\begin{split}
& \theta\triangleq \begin{pmatrix}
      \widetilde{X}   \\
     \widetilde{ \Lambda}_1
\end{pmatrix},~~e_k \triangleq \begin{pmatrix}
      e_{1,k}+e_{2,k}    \\
        ( \mathcal{V}_1^T\otimes \mathbf{I}_m)  e_{3,k}  
\end{pmatrix} ,    \\
& g(\theta)\triangleq-\begin{pmatrix}
      & g_1(\theta)   \\
      & g_2(\theta) 
\end{pmatrix} = -\begin{pmatrix}
     \ \nabla \widetilde{f}(\widetilde{X}+X^*) - \nabla \widetilde{f}(X^*)  +  ( \mathcal{ \bar{L}}  \otimes \mathbf{I}_m)  \widetilde{X}+(\mathcal{V}_1\mathcal{S} \otimes \mathbf{I}_m) \widetilde{\Lambda}_{1}    \\
       -(\mathcal{S} \mathcal{V}_1^T \otimes \mathbf{I}_m) \widetilde{X} \end{pmatrix} .
\end{split}
\end{equation}
Then we can rewrite \eqref{dual2} \eqref{primal2} as
$$\theta_{k+1}=\theta_k+\gamma_kY_k,$$
where  $Y_k=g(\theta_k)+e_k$.

 We  want to apply Lemma \ref{Stable1} i). For this, we have to   validate  conditions C0-C3.

\textbf{Step 1:}  We first  show   C0.
  By Assumption \ref{ass-set}-b, from   \cite[Theorem 2.1.5]{Nesterov} it follows that 
  \begin{equation}\label{equc}
\langle  x-y , \nabla f_i(x) - \nabla f_i(y)\rangle  \geq  \frac{1}{L_f } \| \nabla f_i(x)- \nabla f_i(y) \|^2 ~~\forall x,y \in \mathbb{R}^m.
\end{equation}
 
Set $V_1(\theta) \triangleq \frac{1}{2}(\|\widetilde{X}\|^2+\|\widetilde{\Lambda}_1\|^2).$ 
Then by \eqref{gtheta} \eqref{equc} we obtain  \begin{equation}\label{Lyv1}
\begin{split}
&\nabla V_1(\theta)^T g(\theta)\\&=- \widetilde{X}^T\big(  \ \nabla \widetilde{f}(\widetilde{X}+X^*) - \nabla \widetilde{f}(X^*)  +  ( \mathcal{ \bar{L}}  \otimes \mathbf{I}_m)  \widetilde{X}+(\mathcal{V}_1\mathcal{S} \otimes \mathbf{I}_m)\widetilde{\Lambda}_{1} \big)+ \widetilde{\Lambda}_1^T  (\mathcal{S} \mathcal{V}_1^T \otimes \mathbf{I}_m) \widetilde{X} \\&=- \widetilde{X}^T ( \mathcal{ \bar{L}}  \otimes \mathbf{I}_m)  \widetilde{X}-\alpha\widetilde{X}^T\big(  \ \nabla \widetilde{f}(\widetilde{X}+X^*) - \nabla \widetilde{f}(X^*)\big)\\&
\leq - \widetilde{X}^T ( \mathcal{ \bar{L}}  \otimes \mathbf{I}_m)  \widetilde{X}-  \frac{1}{L_f } \|  \ \nabla \widetilde{f}(\widetilde{X}+X^*) - \nabla \widetilde{f}(X^*)\|^2\\&
\leq - \widetilde{X}^T ( \mathcal{ \bar{L}}  \otimes \mathbf{I}_m)  \widetilde{X}.
\end{split}
\end{equation}

Set $V_2(\theta) \triangleq \widetilde{f}(\widetilde{X}+X^*)-\widetilde{f}(X^*)- \widetilde{X}^T \nabla \widetilde{f}(X^*) +
 \frac{1}{2}\big( \widetilde{X}^T( \mathcal{ \bar{L}}  \otimes \mathbf{I}_m)  \widetilde{X} \big)+   \widetilde{X}^T(\mathcal{V}_1\mathcal{S} \otimes \mathbf{I}_m)\widetilde{\Lambda}_{1}.$  Then
 \begin{equation}\label{Lyv2}
\begin{split}
&\nabla V_2(\theta)^T g(\theta)=- \|g_1(\theta)\|^2+ \|  (\mathcal{S} \mathcal{V}_1^T \otimes \mathbf{I}_m) \widetilde{X}\|^2. 
\end{split}
\end{equation}
By \eqref{right} we have
$$ (\mathcal{V}^T\mathcal{\bar{L}} \otimes \mathbf{I}_m) \widetilde{X}  = col\{
          (\mathcal{S} \mathcal{V}_1^T \otimes \mathbf{I}_m) \widetilde{X},   \mathbf{0}\}.$$
Then by $\mathcal{V}\mathcal{V}^T=\mathbf{I}_n$ and the properties of Kronecker products \eqref{pkrp} we derive
  $$ \|  (\mathcal{S} \mathcal{V}_1^T \otimes \mathbf{I}_m) \widetilde{X}\|^2 =\| (\mathcal{V}^T\mathcal{\bar{L}} \otimes \mathbf{I}_m) \widetilde{X} \|^2= \widetilde{X} ^T(\mathcal{\bar{L}}^2\otimes \mathbf{I}_m) \widetilde{X} .$$
Hence by \eqref{Lyv2} we derive
 \begin{equation}\label{Lyv3}
\begin{split}
&\nabla V_2(\theta)^T g(\theta)=- \|g_1(\theta)\|^2+ \widetilde{X} ^T(\mathcal{\bar{L}}^2\otimes \mathbf{I}_m) \widetilde{X} . 
\end{split}
\end{equation}

Set $V(\theta)\triangleq V_1(\theta)+\alpha V_2(\theta) $ with  $0< \alpha<\frac{1}{ \kappa^*  }$, where $\kappa^*= \max\limits_{i=2,\cdots,n}\kappa_i $.  
Then by \eqref{Lyv1} and \eqref{Lyv3} we derive
 \begin{equation}\label{Lyv4}
\begin{split}
&\nabla V(\theta)^T g(\theta)=-  \widetilde{X} ^T\big((\mathcal{\bar{L}}-\alpha\mathcal{\bar{L}}^2)\otimes \mathbf{I}_m\big) \widetilde{X} -\alpha \|g_1(\theta)\|^2  . 
\end{split}
\end{equation}
Since   $\mathcal{V}^T\mathcal{L} \mathcal{V}=diag\{0,\kappa_2,  \cdots, \kappa_n\},$ we have $\mathcal{V}^T\mathcal{\bar{L}}^2 \mathcal{V}=diag\{0,\kappa_2^2, \cdots, \kappa_n^2\}  $.
  Then   all  possible distinct  eigenvalues of    $\mathcal{\bar{L}}-\alpha\mathcal{\bar{L}}^2$ are $0, $ and $ \kappa_i-\alpha\kappa_i^2, i=2,\cdots, n$. By $ 0< \alpha \leq \frac{1}{\kappa_n}  $ we derive   $\alpha\kappa_i \leq 1~~ \forall i=1,\cdots, n$,  and hence
  $  \kappa_i-\alpha\kappa_i^2=  \kappa_i(1- \alpha\kappa_i)\geq 0~ ~\forall i=1,\cdots, n.$
 Thus  for  any $\alpha$ with  $0< \alpha<\frac{1}{ \kappa^*  }$,     the   matrix  $\mathcal{L} - \alpha\mathcal{L}^2$  is positive semi-definite.  Then by \eqref{Lyv4} we have
 \begin{equation}\label{Lyv5}
\nabla V(\theta)^T g(\theta) \leq 0~~\forall \theta \in \mathbb{R}^{(2n-1)  m}.
\end{equation}
 The equality holds if and only if $ \widetilde{X} ^T\big((\mathcal{\bar{L}}-\alpha\mathcal{\bar{L}}^2)\otimes \mathbf{I}_m\big) \widetilde{X} =0,~g_1(\theta)=\mathbf{0}.$
 
 Since the   matrix  $\mathcal{L} - \alpha\mathcal{L}^2$  is positive semi-definite, the equality
 $ \widetilde{X} ^T\big((\mathcal{\bar{L}}-\alpha\mathcal{\bar{L}}^2)\otimes \mathbf{I}_m\big) \widetilde{X} =0$ implies that  $\widetilde{X}=\mathbf{1} \otimes \tilde{x}$. Then  by multiplying both sides of  
  \begin{equation}\label{equ1} \nabla \widetilde{f}(\widetilde{X}+X^*) - \nabla \widetilde{f}(X^*)  +  ( \mathcal{ \bar{L}}  \otimes \mathbf{I}_m)  \widetilde{X}+(\mathcal{V}_1\mathcal{S} \otimes \mathbf{I}_m)  \widetilde{\Lambda}_{1} =\mathbf{0} 
\end{equation}  from  left with   $\mathbf{1}^T\otimes \mathbf{I}_m $,  from $\mathbf{1}^T\mathcal{V}_1=0$ and $\mathbf{1}^T\mathcal{\bar{L}}=0$ by  \eqref{pkrp} it follows that   $$\nabla f(x^*+\tilde{x})-\nabla f(x^*)=\mathbf{0}.$$
  Since $f(\cdot)$ is strictly convex with  $x^*$  being the  unique optimal solution, 
  $\nabla f(x^*+\tilde{x})= \nabla f(x^*)=\mathbf{0},$  and hence, $\tilde{x}=\mathbf{0}$. Then from \eqref{equ1} we see $(\mathcal{V}_1\mathcal{S} \otimes \mathbf{I}_m ) \widetilde{\Lambda}_{1} =\mathbf{0} $,   and hence $(\mathcal{V}_1^T\mathcal{V}_1\mathcal{S} \otimes \mathbf{I}_m ) \widetilde{\Lambda}_{1} =\mathbf{0} $. 
 By noticing that $\mathcal{V}_1^T\mathcal{V}_1=\mathbf{I}_{n-1}$ and   $S$  is  a diagonal matrix with positive diagonal entries,   we obtain $\widetilde{\Lambda}_{1} =\mathbf{0} $. Consequently, $\nabla V(\theta)^T g(\theta) =0$ only if  $  \theta=\mathbf{0}$.  Therefore,  by  \eqref{Lyv5} we derive C0.

   \textbf{Step 2: } We now verify C1.   We use    Lemma  \ref{Stable0}    to prove    $ \lim\limits_{k\rightarrow \infty} \theta_k =\mathbf{0}   ~a.s.$

Note that $X_k,\Lambda_k$ are bounded with probability one by Theorem \ref{thm2}.  
Then by the definition of $\theta_k $ we know that C1' holds. 
From Lemma \ref{lemma1} and Assumption \ref{ass-stepsize1}, by the convergence theorem for mds \cite{YSChow}  $$\sum_{k=1}^{\infty} \gamma_k e_{k}<\infty~a.s.,$$  
and hence C2' holds.
By the definition of $g(\theta)$ it is seen that C3' holds.  
Since it has already been  proven in Step 1 that C0 holds,  by Lemma  \ref{Stable0}    we obtain C1.

  \textbf{Step 3:}  We now verify  C2. 
Define $\varepsilon_k=e_kI_{[\|\theta_k \|\leq \epsilon]},~ \nu_k=e_kI_{[\|\theta_k \|> \epsilon]}$, where $\epsilon>0$ is a   constant. 

By noting that  $ \lim\limits_{k\rightarrow \infty} \theta_k =\mathbf{0}   ~a.s.$, 
there exists $k_0$ possibly depending on samples  such that 
\begin{equation}\label{unib}
\|\theta_k \| \leq \epsilon~~\forall k \geq k_0~~a.s.
\end{equation}Thus,  
$\nu_k=\mathbf{0}~\forall k \geq k_0~a.s.,$
and hence \eqref{nuk} holds. 

  Since $\theta_k$ is adapted to  $ \mathcal{F}_{k-1}$, from   \eqref{noise1}, \eqref{noise2}, \eqref{noise3}  
  and by $e_k$ defined in \eqref{gtheta} we derive  
\begin{equation}\label{csm}
E[\varepsilon_k | \mathcal{F}_{k-1}]=E[e_k | \mathcal{F}_{k-1}]I_{[\|\theta_k \|\leq \epsilon]}=\mathbf{0}~a.s.
\end{equation}

Since   $(\mathcal{ L} _k \otimes \mathbf{I}_m ) X^*=\mathbf{0}$,  by \eqref{def-e1} we derive 
\begin{equation} 
\begin{split}
  e_{1,k} &= \big( (\mathcal{\bar{L}}-\mathcal{ L} _k)\otimes \mathbf{I}_m  \big)
  \big( \Lambda_k-\Lambda^*  \big)    +\big( (\mathcal{\bar{L}}-\mathcal{ L} _k)\otimes \mathbf{I}_m  \big)\widetilde{X}_k  +\big( (\mathcal{\bar{L}}-\mathcal{ L} _k )\otimes \mathbf{I}_m \big) \Lambda^*.\nonumber
\end{split}
\end{equation}
By noticing that $ \mathcal{V}_1\mathcal{V}_1^T =\mathbf{I}_n-\frac{\mathbf{1}\mathbf{1}^T}{n}$ 
and $\mathcal{\bar{L}}\mathbf{1}= \mathcal{L}_k \mathbf{1}=\mathbf{0} $, we derive
\begin{equation}\label{mult1}
\begin{split}
    \mathcal{\bar{L}}-\mathcal{ L} _k  =( \mathcal{\bar{L}}-\mathcal{ L} _k )\mathcal{V}_1\mathcal{V}_1^T. 
\end{split}
\end{equation}
Then by $\widetilde{\Lambda}_{1,k}=(\mathcal{V}_1^T \otimes \mathbf{I}_m)(\Lambda_k-\Lambda^*)$ we see that 
\begin{equation}\label{lambda1}
\begin{split}
  &  (\mathcal{\bar{L}}-\mathcal{ L} _k )\otimes \mathbf{I}_m 
  ( \Lambda_k-\Lambda^* ) = ( \mathcal{\bar{L}}-\mathcal{ L} _k )\mathcal{V}_1 \otimes \mathbf{I}_m \widetilde{\Lambda}_{1,k}. 
\end{split}
\end{equation}

Then by multiplying both sides of the first equality in \eqref{optimalcd} with $\mathcal{V}^T \otimes \mathbf{I}_m $ from  left, and  by  \eqref{pkrp} \eqref{right}  we obtain
\begin{equation}
\label{ }
\begin{split}
  & -(\mathcal{V}^T \otimes \mathbf{I}_m)\nabla \widetilde{f}(X^*) = (\mathcal{V}^T\otimes \mathbf{I}_m)  ( \mathcal{ \bar{L}}  \otimes \mathbf{I}_m) \Lambda^* = \big(\mathcal{V}^T  \mathcal{ \bar{L}}  \otimes \mathbf{I}_m\big) \Lambda^* =col\{
    ( \mathcal{S} \mathcal{V}_1^T \otimes \mathbf{I}_m)\Lambda^*,
 \mathbf{0 }  \}.
   \nonumber
\end{split}
\end{equation}
So, $ (\mathcal{S} \mathcal{V}_1^T\otimes \mathbf{I}_m)  \Lambda^* =-(\mathcal{V}_1^T\otimes \mathbf{I}_m)\nabla \widetilde{f}(X^*).$ 
Then by \eqref{mult1} we obtain
\begin{equation}\label{mult2}
\begin{split}
 & \big( (\mathcal{\bar{L}}-\mathcal{ L} _k)\otimes \mathbf{I}_m  \big) \Lambda^*  =\big(( \mathcal{\bar{L}}-\mathcal{ L} _k )\mathcal{V}_1\mathcal{V}_1^T\otimes \mathbf{I}_m \big ) \Lambda^*
  \\&=\big(( \mathcal{\bar{L}}-\mathcal{ L} _k )\mathcal{V}_1S^{-1} S   \mathcal{V}_1^T \otimes \mathbf{I}_m\big) \Lambda^*  =\big((  \mathcal{ L} _k -   \mathcal{ \bar{L}} ) \mathcal{V}_1S^{-1} \mathcal{V}_1^T\otimes \mathbf{I}_m\big)\nabla \widetilde{f}(X^*).  \nonumber
 \end{split}
\end{equation}
Thus, by \eqref{lambda1} we obtain
\begin{equation}\label{e1}
\begin{split}
  e_{1,k}&= \big( ( \mathcal{\bar{L}}-\mathcal{ L} _k )\mathcal{V}_1 \otimes \mathbf{I}_m)\big) \widetilde{\Lambda}_{1,k}    +\big( (\mathcal{\bar{L}}-\mathcal{ L} _k )\otimes \mathbf{I}_m \big)\widetilde{X}_k   \\&+\big((  \mathcal{ L} _k -   \mathcal{ \bar{L}} ) \mathcal{V}_1S^{-1} \mathcal{V}_1^T\otimes \mathbf{I}_m\big)\nabla \widetilde{f}(X^*). 
\end{split}
\end{equation}

Note that   $X_k   $ and  $ \Lambda_k$ are adapted to $\mathcal{F}_{k-1}  $, and   that  $\mathcal{L}_k$
is  independent of  $\mathcal{F}_{k-1}$  by Assumption \ref{ass-graph}-d. Then by \eqref{e1}  we  derive
\begin{equation}\label{noise-cond1}
\begin{split}
E[ \| e_{1,k}\|^2 | \mathcal{F}_{k-1}]  \leq    3C_{01}( \| \widetilde{ X}_k\|^2+ \|\widetilde{\Lambda}_{1,k}\|^2+C_{04} )~~a.s.,         \end{split}\end{equation}
where $C_{04}=\|(\mathcal{V}_1S^{-1} \mathcal{V}_1^T \otimes \mathbf{I}_m)\nabla \widetilde{f}(X^*)\|^2.$  
Since  $X_k=X^*+\widetilde{X}_k$,  by   \eqref{noise2}\eqref{noise3} we derive \begin{align}
& E[\| e_{2,k}\|^{2}  | \mathcal{F}_{k-1}] \leq C_{02} +6c_v (\|X^*\|^2+ \| \widetilde{ X}_k\|^2)~~a.s.,  \label{noise-cond2}   \\&
E[\| e_{3,k}\|^{2}  | \mathcal{F}_{k-1}]\leq  2C_{01}(\|X^*\|^2+ \| \widetilde{ X}_k\|^2) + C_{03}~~a.s. \label{noise-cond3}
\end{align}

Since $\theta_k$ is adapted to $\mathcal{F}_{k-1}$, by \eqref{noise-cond1} \eqref{noise-cond2} \eqref{noise-cond3} we know that there exists a constant $K>0$ such that    \begin{equation}\label{condnoise}
E[\| \varepsilon_k\|^{2}  | \mathcal{F}_{k-1} ]=E[\| e_{k}\|^{2}| \mathcal{F}_{k-1} ] I_{[ \| \theta_k  \|\leq \epsilon]} \leq  K~~\forall k\geq 1 ~~a.s.
\end{equation}
 Consequently,  by \eqref{csm} and \eqref{condnoise} we  know  that \eqref{A-1} holds.
 
By the Chebyshev's  inequality  from \eqref{condnoise} we have 
$$\mathbb{P}(\|\varepsilon_k \| > a) \leq \frac{E[\| \varepsilon_k \|^2]}{a^2} \leq \frac{K}{a^2}~~\forall k \geq 1.$$
Then by the  Schwarz inequality from \eqref{condnoise} we derive 
\begin{equation} 
 \begin{split}   E[\|\varepsilon_k I_{[\|\varepsilon_k \| > a]} \|] & \leq   (E[\|\varepsilon_k  \|^2] )^{\frac{1}{2}}(E[I_{[\varepsilon_k \| > a]}^2] )^{\frac{1}{2}}  \leq \sqrt{K} \sqrt{\mathbb{P}(\|\varepsilon_k \| > a)} \leq  \frac{K}{a} ~~\forall k \geq 1.
    \nonumber
\end{split}
\end{equation}
Therefore,
$\lim\limits_{a \rightarrow \infty} \sup_k E[\|\varepsilon_k I_{[\|\varepsilon_k \| > a]} \|]=0,$ and hence \eqref{A-3} holds.
  
  Note that
 \begin{equation}\label{defe}
\begin{split}
 e_ke_k^T&=\begin{pmatrix}
      e_{1,k}+e_{2,k}    \\
         (\mathcal{V}_1^T  \otimes \mathbf{I}_m ) e_{3,k}  
\end{pmatrix}  (    e_{1,k}^T+e_{2,k} ^T,    e_{3,k}^T   (\mathcal{V}_1   \otimes \mathbf{I}_m ) \\
&= \begin{pmatrix}
      &    (e_{1,k}+e_{2,k} )  (e_{1,k}+e_{2,k} ) ^T  & (e_{1,k}+e_{2,k} ) e_{3,k}^T( \mathcal{V}_1 \otimes \mathbf{I}_m ) \\
      &     (\mathcal{V}_1^T \otimes \mathbf{I}_m ) e_{3,k}   (e_{1,k}+e_{2,k} )^T &    (\mathcal{V}_1^T \otimes \mathbf{I}_m ) e_{3,k}  e_{3,k}^T ( \mathcal{V}_1 \otimes \mathbf{I}_m )\end{pmatrix},
\end{split}
\end{equation}
and  that $ \lim\limits_{k\rightarrow \infty} \widetilde{ X}_k =\mathbf{0}~ a.s. , \lim\limits_{k\rightarrow \infty} \widetilde{  \Lambda}_{1,k} =\mathbf{0} ~ a.s. $, and  $\widetilde{ X}_k,\widetilde{  \Lambda}_{1,k}$ are adapted to $\mathcal{F}_{k-1}$. Then by  Assumption \ref{ass-graph}-d, from \eqref{unib}  \eqref{e1}  and the definition of $S_1$  given in \eqref{symbol} we derive 
 \begin{equation}\label{cond-e1}
\begin{split}
 &E[ e_{1,k}  e_{1,k}^TI_{[ \| \theta_k  \|\leq \epsilon]} | \mathcal{F}_{k-1}] \xlongrightarrow [k \rightarrow \infty]{}  S_1~a.s . \end{split}
\end{equation}

By Assumptions \ref{ass-noise2}-a and \ref{ass-independency}-b we derive  \begin{equation} \label{ind1}
\begin{split}
 E[\omega_{i_1j_1,k}\omega_{i_2j_2,k}^T|\mathcal{F}_{k-1}']&=E[ \omega_{i_1j_1,k} |\mathcal{F}_{k-1}']E[ \omega_{i_2j_2,k}^T|\mathcal{F}_{k-1}']=\mathbf{0}~~\forall (i_1,j_1) \neq (i_2,j_2). \end{split}
\end{equation}
Thus,    noticing that $a_{ij,k}~\forall i,j \in\mathcal{V}$ are adapted to $\mathcal{F}_{k-1}'$ we obtain
 \begin{equation}\label{ }
\begin{split}
 &E[  \omega_{i_1,k}\omega_{i_2,k}^T | \mathcal{F}_{k-1}']  =  \sum_{j_1,j_2=1}^n a_{i_1j_1,k}a_{i_2j_2,k} E[\omega_{i_1j_1,k}  \omega_{i_2j_2,k}^T| \mathcal{F}_{k-1}'] 
  =\mathbf{0}~~\forall i_1 \neq i_2.  \nonumber \end{split}
\end{equation}
Then   by $\mathcal{F}_k \subset\mathcal{F}_k'$    we have 
 \begin{equation}\label{corre-cond-w}
\begin{split}
 &E[  \omega_{i_1,k}\omega_{i_2,k}^T | \mathcal{F}_{k-1}]= E\big[ E[  \omega_{i_1,k}\omega_{i_2,k}^T | \mathcal{F}_{k-1}'] \big| \mathcal{F}_{k-1} \big]
  =\mathbf{0}~~\forall i_1 \neq i_2. \end{split}
\end{equation}

By \eqref{ind1} and Assumption \ref{ass-noise2}-a  we obtain
   \begin{equation}\label{ }
\begin{split}
 E[  \omega_{i,k}\omega_{i,k}^T | \mathcal{F}_{k-1}'] &=  \sum_{j_1,j_2=1}^n a_{ij_1,k}a_{ij_2,k} E[\omega_{ij_1,k}  \omega_{ij_2,k}^T| \mathcal{F}_{k-1}']  = \sum_{j=1}^n a_{ij,k}^2 E[\omega_{ij,k}  \omega_{ij,k}^T| \mathcal{F}_{k-1}'] = \sum_{j=1}^n a_{ij,k}^2 R_{\omega,ij}. \nonumber \end{split}
\end{equation}
Then by  Assumptions \ref{ass-graph}-c and \ref{ass-graph}-d, from $\mathcal{F}_k \subset\mathcal{F}_k'$ it follows that 
 \begin{equation}\label{cond-w}
\begin{split}
 E[  \omega_{i,k}\omega_{i,k}^T | \mathcal{F}_{k-1}]&= E\big[ E[  \omega_{i,k}\omega_{i,k}^T | \mathcal{F}_{k-1}'] \big| \mathcal{F}_{k-1} \big] \\&=E[  \sum_{j=1}^n a_{ij,k}^2 R_{\omega,ij} | \mathcal{F}_{k-1}]
  = \sum_{j=1}^n E[  a_{ij,k}^2 ]R_{\omega,ij}= R_{\omega,i} ~~a.s.\end{split}
\end{equation} 

Similarly,  
\begin{equation}\label{ }
\begin{split}
 &E[  \zeta_{i_1,k}\zeta_{i_2,k}^T | \mathcal{F}_{k-1}]=  \mathbf{0}~~\forall i_1 \neq i_2~~a.s, 
~~E[  \zeta_{i,k}\zeta_{i,k}^T | \mathcal{F}_{k-1}]=    R_{\zeta,i}~~a.s. \nonumber  \end{split}
\end{equation}
From here,  by \eqref{corre-cond-w} and   \eqref{cond-w} we conclude that 
 \begin{equation}\label{cond-wz}
 \begin{split}
 &E[ \omega_k \omega_k^T | \mathcal{F}_{k-1}] =   diag \big \{R_{\omega,1} ,\cdots,R_{\omega,N}  \big\}  =R_{\omega}~~a.s.,  \\& E[\zeta_k \zeta_k^T | \mathcal{F}_{k-1}] =   diag \big \{R_{\zeta,1} ,\cdots,R_{\zeta,N}  \big\}=R_{\zeta}~~a.s. 
 \end{split}
 \end{equation}
By Assumptions \ref{ass-noise2} and  \ref{ass-independency}-c,  similar to  \eqref{cond-wz}  we 
can show that 
\begin{equation}\label{cond-wzv}
\begin{split}
 &E[ \omega_k\zeta_k^T | \mathcal{F}_{k-1}] =   \mathbf{0}~~a.s., ~~ E[ \omega_kv_k^T | \mathcal{F}_{k-1}] =   \mathbf{0}~~a.s., ~~E[ \zeta_k v_k^T | \mathcal{F}_{k-1}] =   \mathbf{0}~~a.s.  \end{split}
 \end{equation}

From  Assumptions \ref{ass-noise2}-c and  \ref{ass-independency}-a it follows that
  $$E[v_{i,k}v_{j,k}^T|\mathcal{F}_{k-1}']=E[v_{i,k} |\mathcal{F}_{k-1}']E[ v_{j,k}^T|\mathcal{F}_{k-1}']=\mathbf{0}~~\forall i \neq j,$$ 
and hence  by Assumption \ref{ass-noise2}-c we obtain 
 \begin{equation}\label{ }
\begin{split}
 E[ v_k v_k^T | \mathcal{F}_{k-1}'] &=   diag \big \{E[ v_{1,k} v_{1,k}^T | \mathcal{F}_{k-1}'] ,\cdots,E[ v_{n,k} v_{n,k}^T | \mathcal{F}_{k-1}'] \big) \\&
  \xlongrightarrow [k \rightarrow \infty]{}   diag \big \{R_{v,1},\cdots,R_{v,n} \big)  = R_v~a.s. \nonumber
   \end{split}
\end{equation}
Noting that  $v_{k}$ and $\mathcal{L}_k$ are conditionally independent  given $\mathcal{F}_{k-1}$ by Assumption   \ref{ass-independency}-d, by   \cite[Corollary 7.3.2]{YSChow}  we have 
$$E[ v_k v_k^T | \mathcal{F}_{k-1}'] =E[ v_k v_k^T | \mathcal{F}_{k-1},\mathcal{L}_k]=E[ v_k v_k^T | \mathcal{F}_{k-1}] \xlongrightarrow [k \rightarrow \infty]{}     R_v~~a.s.$$
 For  $e_{2,k} $  defined by  \eqref{def-e2} ,  by \eqref{cond-wz}  and \eqref{cond-wzv} we derive
 \begin{equation}\label{ }
\begin{split}
 E[ e_{2,k}  e_{2,k}^T  | \mathcal{F}_{k-1}]&= E[ v_kv_k^T +\omega_k\omega_k^T+\zeta_k\zeta_k^T | \mathcal{F}_{k-1}] \xlongrightarrow [k \rightarrow \infty]{}   R_v+R_{\omega}+R_{\zeta}= S_2~~a.s. \nonumber  \end{split}
\end{equation}
Thus by noticing that $\theta_k$ is adapted to  $\mathcal{F}_{k-1}' $, from \eqref{unib} we derive   
 \begin{equation}\label{cond-e2}
\begin{split}
 &E[ e_{2,k}  e_{2,k}^TI_{[  \| \theta_k  \|\leq \epsilon]} | \mathcal{F}_{k-1}]  =E[ e_{2,k}  e_{2,k}^T | \mathcal{F}_{k-1}]   I_{[  \| \theta_k  \|\leq \epsilon]}\xlongrightarrow [k \rightarrow \infty]{}   S_2~~a.s. \end{split}
\end{equation}

    Since $\mathcal{F}_{k} \subset \mathcal{F}_{k}'$ and $e_{1,k},\theta_k$  are adapted to $\mathcal{F}_{k-1}'$, 
   by \eqref{equ3}  we obtain   \begin{equation} \begin{split} &E [ e_{1,k} e_{2,k}^TI_{[ \| \theta_k  \|\leq \epsilon]} | \mathcal{F}_{k-1}]= E \big[  E [  e_{1,k} e_{2,k} ^TI_{[ \| \theta_k  \|\leq \epsilon]}  | \mathcal{F}_{k-1}'] \big | \mathcal{F}_{k-1} \big ] \\&= E \big[  e_{1,k} I_{[ \| \theta_k  \|\leq \epsilon]} E [  e_{2,k}^T  | \mathcal{F}_{k-1}'] \big | \mathcal{F}_{k-1} \big ]=\mathbf{0}~~a.s. ,\nonumber
\end{split}
\end{equation} 
which incorporating with  \eqref{cond-e1}  \eqref{cond-e2} yields  \begin{equation}\label{cov12}
\begin{split}
 &E[ (e_{1,k}+e_{2,k})( e_{1,k}+  e_{2,k})^TI_{[ \| \theta_k  \|\leq \epsilon]} | \mathcal{F}_{k-1}]
    \xlongrightarrow [k \rightarrow \infty]{} S_1+S_2~~a.s. \end{split}
\end{equation}

By \eqref{equ0} we see that
\begin{equation}\label{mdsw}
 E[\omega_{k} | \mathcal{F}_{k-1}'] =\mathbf{0} .
\end{equation}Hence,  noticing that  $\mathcal{F}_{k} \subset \mathcal{F}_{k}'$ and that  $e_{1,k}I_{[ \| \theta_k  \|\leq \epsilon]}$ is adapted to $\mathcal{F}_{k-1}'$,  we obtain 
     \begin{equation}\label{cond1w} 
      \begin{split} 
      E [ e_{1,k} I_{[ \| \theta_k  \|\leq \epsilon]} \omega_k^T| \mathcal{F}_{k-1}] &=E\big[E [ e_{1,k} I_{[ \| \theta_k  \|\leq \epsilon]} \omega_k^T| \mathcal{F}_{k-1}'] \big | \mathcal{F}_{k-1}\big]\\&
      =E\big[ e_{1,k} I_{[ \| \theta_k  \|\leq \epsilon]} E [  \omega_k^T| \mathcal{F}_{k-1}'] \big | \mathcal{F}_{k-1}\big]=\mathbf{0}~~a.s.
\end{split}
\end{equation} 

Note  that  
\begin{equation}\label{tildex}
\begin{split}
 ( \mathcal{L}_k -\mathcal{\bar{L}}) \otimes \mathbf{I}_mX_k =( \mathcal{L}_k -\mathcal{\bar{L}})\otimes \mathbf{I}_m\widetilde{X}_k,
\end{split}
\end{equation}  and that   $\widetilde{ X}_k $ and $ \widetilde{  \Lambda}_{1,k}$ are adapted to $\mathcal{F}_{k-1}$.
Then  from  \eqref{def-e3} \eqref{e1}   \eqref{cond1w}, by Assumption \ref{ass-graph}-c and \ref{ass-graph}-d we derive
 \begin{equation}\label{cond23}
\begin{split}
 &E[ e_{1,k}  e_{3,k}^TI_{[ \| \theta_k  \|\leq \epsilon]} | \mathcal{F}_{k-1}]  =E[ e_{1,k} \big( ( \mathcal{L}_k -\mathcal{\bar{L}})\otimes \mathbf{I}_m\widetilde{X}_k \big)^TI_{[ \| \theta_k  \|\leq \epsilon]} | \mathcal{F}_{k-1}] \xlongrightarrow [k \rightarrow \infty]{}  \mathbf{0}~~a.s.,   \end{split}
\end{equation}
where the limit takes place because   $ \lim\limits_{k\rightarrow \infty}\widetilde{ X}_k =\mathbf{0}~a.s.$, and $    \lim\limits_{k\rightarrow \infty} \widetilde{  \Lambda}_{1,k}  =\mathbf{0}   ~a.s.$

Since    $\mathcal{L}_k, X_k $ are adapted to $\mathcal{F}_{k-1}'$,   by \eqref{def-e3}\eqref{equ3}     we derive
   \begin{equation}  \begin{split}  E [ e_{2,k} e_{3,k}^T| \mathcal{F}_{k-1}'] & = E [  e_{2,k}  \big(( \mathcal{L}_k -\mathcal{\bar{L}})  \otimes \mathbf{I}_m  X_k  \big)^T  | \mathcal{F}_{k-1}'  ]- E [  e_{2,k}  \omega_k^T  | \mathcal{F}_{k-1}'  ] \\&
   = E [  e_{2,k}   | \mathcal{F}_{k-1}'  ]  \big(( \mathcal{L}_k -\mathcal{\bar{L}}) \otimes \mathbf{I}_m  X_k  \big)^T- E [  e_{2,k}  \omega_k^T  | \mathcal{F}_{k-1}'  ]  = - E [  e_{2,k}  \omega_k^T  | \mathcal{F}_{k-1}'  ] . \nonumber
\end{split}
\end{equation} 
Then by $\mathcal{F}_{k} \subset \mathcal{F}_{k}'$ we conclude that
  \begin{equation}  \begin{split} E [ e_{2,k} e_{3,k}^T| \mathcal{F}_{k-1}]&= E \big[  E [  e_{2,k} e_{3,k} ^T  | \mathcal{F}_{k-1}'] \big | \mathcal{F}_{k-1} \big ] \\&= -E \big[ E [  e_{2,k}  \omega_k^T  | \mathcal{F}_{k-1}'  ]  | \mathcal{F}_{k-1} \big ]=-E [  e_{2,k}  \omega_k^T  | \mathcal{F}_{k-1}]~~a.s. \nonumber
\end{split}
\end{equation} 
Noticing $ e_{2,k}$  defined by \eqref{def-e2},  by \eqref{cond-wz} and \eqref{cond-wzv}  we derive 
$$E [ e_{2,k} e_{3,k}^T| \mathcal{F}_{k-1}]=-E [  \omega_k  \omega_k^T  | \mathcal{F}_{k-1}]=-R_{\omega}~~a.s.$$
Since  $\theta_k$ is adapted to $\mathcal{F}_{k-1}$,  by \eqref{unib} we obtain 
  \begin{equation} \begin{split}
 & \lim_{k \rightarrow \infty}E  [ e_{2,k} e_{3,k}^TI_{[ \| \theta_k  \|\leq \epsilon]} | \mathcal{F}_{k-1}] 
 = \lim_{k \rightarrow \infty}E  [ e_{2,k} e_{3,k}^T | \mathcal{F}_{k-1}]I_{[ \| \theta_k  \|\leq \epsilon]}
  =-R_{\omega} \lim_{k \rightarrow \infty}I_{[ \| \theta_k  \|\leq \epsilon]} =-R_{\omega}~~a.s.,  \nonumber \end{split}
\end{equation} 
which incorporating with \eqref{cond23} yields
 \begin{equation}\label{cov123}
\begin{split}
 &E[ (e_{1,k}+e_{2,k} ) e_{3,k}^T ( \mathcal{V}_1 \otimes \mathbf{I}_m )I_{[ \| \theta_k  \|\leq \epsilon]} | \mathcal{F}_{k-1}]   \xlongrightarrow [k \rightarrow \infty]{}  -R_{\omega} ( \mathcal{V}_1 \otimes \mathbf{I}_m )~~a.s. \end{split}
\end{equation}

Since $\mathcal{L}_k$ and $\widetilde{X}_k$ are adapted to $\mathcal{F}_{k-1}'$, by \eqref{mdsw} and   $\mathcal{F}_{k} \subset \mathcal{F}_{k}'$  we obtain 
     \begin{equation} \label{condwx}
      \begin{split}
      & E [ \omega_k \big(( \mathcal{L}_k -\mathcal{\bar{L}}) \otimes \mathbf{I}_m\widetilde{X}_k \big)^T| \mathcal{F}_{k-1}] \\&=
     E\big[ E [ \omega_k \big(( \mathcal{L}_k -\mathcal{\bar{L}}) \otimes \mathbf{I}_m\widetilde{X}_k \big)^T| \mathcal{F}_{k-1}'] \big| \mathcal{F}_{k-1}\big]\\&=   E\big[ E [ \omega_k | \mathcal{F}_{k-1}']  \big(( \mathcal{L}_k -\mathcal{\bar{L}})\otimes \mathbf{I}_m\widetilde{X}_k \big)^T\big| \mathcal{F}_{k-1}\big]=\mathbf{0}~~a.s.
\end{split}
\end{equation} 
By  definition of $e_{3,k}$ and \eqref{tildex} we see
\begin{equation}\label{cov3}
\begin{split}
   e_{3,k}  e_{3,k}^T &=-( \mathcal{L}_k -\mathcal{\bar{L}})  \otimes \mathbf{I}_m\widetilde{X}_k \omega_k^T-\omega_k \big(( \mathcal{L}_k -\mathcal{\bar{L}})  \otimes \mathbf{I}_m\widetilde{X}_k \big)^T \\ 
   &+\big( ( \mathcal{L}_k -\mathcal{\bar{L}})  \otimes \mathbf{I}_m \big)\widetilde{X}_k \widetilde{X}_k^T \big(( \mathcal{L}_k -\mathcal{\bar{L}})  ^T  \otimes \mathbf{I}_m\big) +\omega_k\omega_k^T.  \end{split}
\end{equation}

Since $ \lim\limits_{k\rightarrow \infty} \widetilde{ X}_k =\mathbf{0}~a.s.,$ and  $\widetilde{ X}_k$ is adapted to $\mathcal{F}_{k-1},$   by Assumption \ref{ass-graph}-c and  \ref{ass-graph}-d we obtain
\begin{equation}\label{ }
\begin{split}
  & \| E\big[  ( \mathcal{L}_k -\mathcal{\bar{L}})  \otimes \mathbf{I}_m  \widetilde{X}_k \widetilde{X}_k^T  ( \mathcal{L}_k -\mathcal{\bar{L}})  ^T  \otimes \mathbf{I}_m \big| \mathcal{F}_{k-1}\big]\| 
    \leq E[\| \mathcal{L}_k -\mathcal{\bar{L}}\|^2] \| \widetilde{ X}_k \|^2\xlongrightarrow [k \rightarrow \infty]{}  0 ~~a.s.\nonumber \end{split}
\end{equation}
Then by   \eqref{cond-wz}  \eqref{condwx}   \eqref{cov3} we obtain
\begin{equation}\label{ }
\begin{split}
 &E[ e_{3,k}  e_{3,k}^TI_{[ \| \theta_k  \|\leq \epsilon]} | \mathcal{F}_{k-1}]  \xlongrightarrow [k \rightarrow \infty]{}  R_{\omega}~~a.s.,\nonumber \end{split}
\end{equation}
which incorporating with \eqref{defe}, \eqref{cov12} and \eqref{cov123} yields 
\begin{equation}\label{ }
\begin{split}
 &E[ e_{k}  e_{k}^TI_{[ \| \theta_k  \|\leq \epsilon]} | \mathcal{F}_{k-1}]  \xlongrightarrow [k \rightarrow \infty]{}   \Sigma_1~~a.s.\nonumber \end{split}
\end{equation}
Hence  by the definition of $ \varepsilon_k $ we obtain
\begin{equation}\label{tidle2}
E[ \varepsilon_k  \varepsilon_k^T  | \mathcal{F}_{k-1}]  \xlongrightarrow [k \rightarrow \infty]{}   \Sigma_1~~a.s.
\end{equation}

By \eqref{condnoise} we derive
$$ E\big[\sup_k E[\| \varepsilon_k\|^{2}  | \mathcal{F}_{k-1} ]\big] \leq K.$$ 
Then by the Lebesgue dominated convergence theorem \cite[Corroloary 4.2.3]{ YSChow} and  by \eqref{tidle2}   we have
\begin{equation}
\begin{split}
\lim_{k \rightarrow \infty}E\big[E[ \varepsilon_k  \varepsilon_k^T  | \mathcal{F}_{k-1}] \big]=E\big[ \lim_{k \rightarrow \infty} E[ \varepsilon_k  \varepsilon_k^T  | \mathcal{F}_{k-1}] \big]= \Sigma_1. \nonumber
\end{split}
\end{equation}
Thus, $\lim\limits_{k \rightarrow \infty}E[ \varepsilon_k  \varepsilon_k^T ]=\lim\limits_{k \rightarrow \infty}E\big[E[ \varepsilon_k  \varepsilon_k^T  | \mathcal{F}_{k-1}] \big]=\Sigma_1,$ and hence \eqref{A-2} holds.
So, C2 has been verified.

Step 4:  It  remains to check   C3. By \eqref{gtheta} and \eqref{hessian} we derive
\begin{equation} 
\begin{split}
   g(\theta)-F\theta &=   \begin{pmatrix}
      \big( ( \mathcal{ \bar{L}}  \otimes \mathbf{I}_m)  +\mathcal{H}\big) \widetilde{X} + \mathcal{V}_1\mathcal{S} \otimes \mathbf{I}_m \widetilde{\Lambda}_1\\
    -\mathcal{S} \mathcal{V}_1^T \otimes \mathbf{I}_m   \widetilde{X}
\end{pmatrix} \\& -\begin{pmatrix}
    \nabla \widetilde{f}(\widetilde{X}+X^*) - \nabla \widetilde{f}(X^*)  +  ( \mathcal{ \bar{L}}  \otimes \mathbf{I}_m)  \widetilde{X}+(\mathcal{V}_1\mathcal{S} \otimes \mathbf{I}_m) \widetilde{\Lambda}_{1}    \\
       -(\mathcal{S} \mathcal{V}_1^T \otimes \mathbf{I}_m) \widetilde{X} \end{pmatrix}
       \\&=- \begin{pmatrix}
      &  \nabla \widetilde{f}(\widetilde{X}+X^*) - \nabla \widetilde{f}(X^*) -  \mathcal{H}  \widetilde{X}   \\
      &  \mathbf{0}
\end{pmatrix} . \nonumber 
\end{split}
\end{equation}
Then by  Assumption \ref{ass-function}-c we obtain  
$$\|  g(\theta)-F\theta\|^2\leq c \|\widetilde{X}\|^2\leq  c \|\theta\|^2,$$
 and hence C3 holds. 
 
 In summary, we have verified C0-C3. Then by   Lemma \ref{Stable1} i) the assertion of the theorem follows. 
\hfill  $\blacksquare$

 \textbf{ Proof of Theorem \ref{thm4}:}    
Since it is  shown in the proof of Theorem \ref{thm3} that C0-C3 hold,   by Lemma \ref{Stable1} ii)  we immediately   derive the assertion.  \hfill  $\blacksquare$

   \section{Numerical Examples} \label{sec:Simulation}

   In this section, we do  simulations for    the   distributed parameter estimation problem considered in \cite{Sayed_TSP_2015}.   We aim at   estimating  the    unknown  $m $-dimensional vector
 $x^*$ based on the data gathered by  $n$ spatially distributed sensors  in the  network.
 Each agent    $i=1,\cdots,n$  at time $k$ has access to  its real scalar  measurement $ d_{i,k}$ 
  given by  the following linear time-varying  model
 \begin{equation}\label{model}
    d_{i,k}=u_{i,k}x^* +\nu_{i,k}, \nonumber
 \end{equation}
  where   $u_{i,k} \in  \mathbb{R}^{1  \times m}$ is the regression vector  accessible to  agent $i$,
  and  $\nu_{i,k}  $   is the local observation noise of   agent $i$.

Assume that $\{u_{i,k}\}$ and $\{\nu_{i,k}\}$ are mutually independent   iid  Gaussian sequences with  distributions $N(\mathbf{0}, R_{u,i})$ and    $N(0,\sigma_{i,\nu}^2)$, respectively.   Besides, we allow some   covariance matrices   nonpositive definite, but  require $\sum_{i=1}^n R_{u,i}$    be positive definite. 
This parameter estimation problem is modeled as solving the following distributed stochastic optimization problem 
 \begin{equation}\label{filter1}
  \min_{x} ~  f(x)=\sum_{i=1}^n f_i(x)\deq E   [   \parallel  d_{i,k}-u_{i,k}x\parallel ^2]. 
  \end{equation}
   So, $f_i(x)=(x-x^*)^T R_{u,i}(x-x^*)+\sigma_{i,\nu},^2$ and $\nabla f_i(x)= R_{u,i}(x-x^*).$
  Therefore,   $x^*$ is the unique optimal solution to  \eqref{filter1} when $\sum_{i=1}^n R_{u,i}$  is positive definite.

Let $m=3$. Set $x^*=(1,2,3)$, and   $$R_{u,1}=\begin{pmatrix}
    1 & 0 & 0 \\
0 &  1&0 \\
      0  & 0&0
\end{pmatrix},~R_{u,2}=\begin{pmatrix}
    0 & 0 & 0 \\
0 &  1&0 \\
      0  & 0&1
\end{pmatrix} ,R_{u,2}=\begin{pmatrix}
   1 & 0 & 0 \\
0 &  0&0 \\
      0  & 0&1
\end{pmatrix},~ \sigma_{i,\nu}=\sqrt{ 0.1}~\forall i \in \mathcal{V}.$$ 
Set  $n=3$  with the underling undirected  graph being  fully connected.
At any  time $k\geq 0$,   with equal probability  $\frac{1}{3}$ for each edge,  we randomly  choose one edge   from the   graph. Set $a_{ij,k}=a_{ji,k}=1$  when  the  edge between $i$ and $j$ is chosen.   For any $i,j\in \mathcal{V}$,
let the communication noises $\{\omega_{ij,k}\}$
and  $\{\zeta_{ij,k}\}$ be mutually independent  iid Gaussian  sequences $N(0,0.1\mathbf{I}_3)$.

   \begin{figure}
       \centering
  \includegraphics[width=5in]{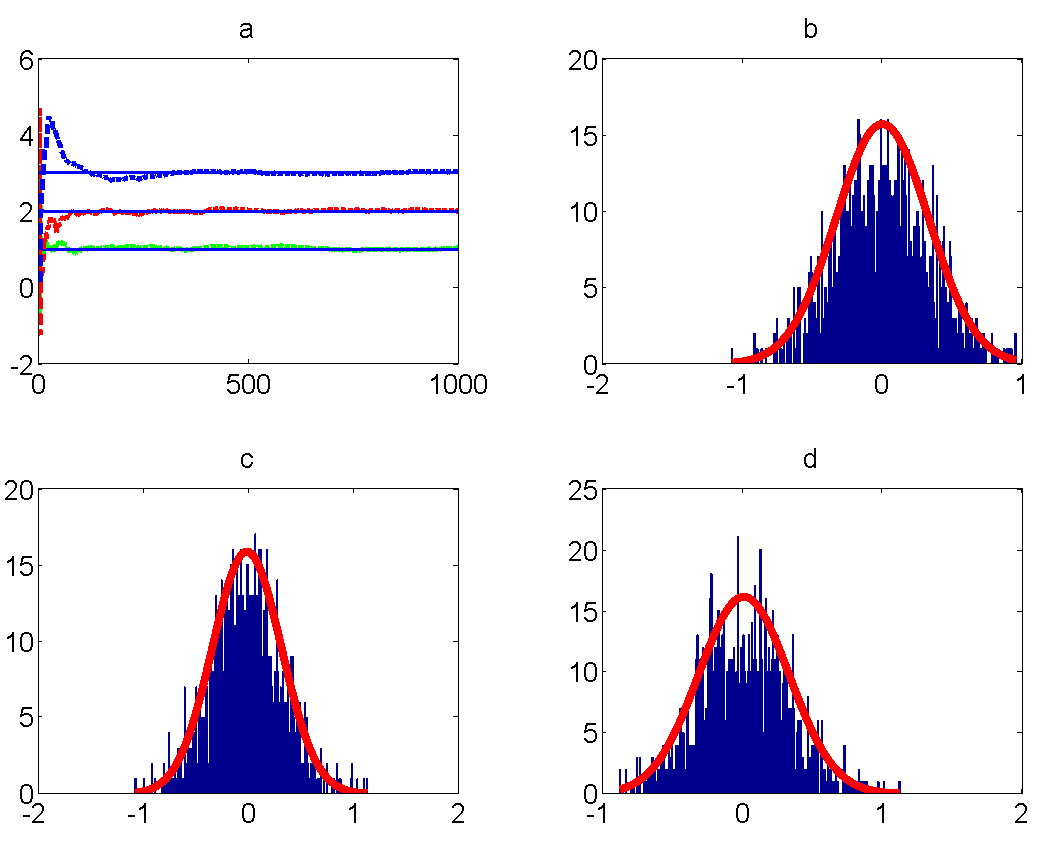}
   \caption{  Estimates $x_{1,k}$,    histographs and    limit distributions for  
    $(x_{1,k}-x^*)/\sqrt{\gamma_k} $ at  $k=1000$} \label{GOne}
\end{figure}

Set $\gamma_k=\frac{1}{k^{0.75}}$.   
By  using     $u_{i,k}$ and $d_{i,k}$ observed at time $k$,  the noisy observation  of the gradient  $\nabla f_i(x_{i,k}) $ is  constructed as $g_{i,k}=u_{i,k}^Tu_{i,k}x_{i,k}-d_{i,k}u_{i,k}^T$.
 Let $\{ x_{i,k}\}$ and  $\{ \lambda_{i,k}\}$   be produced by the algorithm \eqref{algorithm1} with  initial values $x_{i,0}=\mathbf{0},~ \lambda_{i,0}=\mathbf{0}$.  
Since $v_{i,k}=(u_{i,k}^Tu_{i,k}- R_{u,i})(x_{i,k}-x^*)-\nu_{i,k}u_{i,k}^T,$  it is seen that
$\{v_{i,k}\}$ satisfies Assumption \ref{ass-noise}-c
with  $c_v= max\{ E[\|u_{i,k}^Tu_{i,k}- R_{u,i}\|^2], \sigma_{i,v}^2\|R_{u,i}\|\} $.  Then
$\lim\limits_{k \rightarrow \infty} x_{i,k}=x^*$  by Theorem \ref{thm1},
  and hence $\lim\limits_{k \rightarrow \infty} E[v_{i,k} v_{i,k}^T | \mathcal{F}_{k-1}']=\sigma_{i,v}^2 R_{u,i}$. 
As a result, the gradient observation noise for the  distributed parameter estimation problem satisfies 
Assumption \ref{ass-noise2}-c.
 
\begin{figure}
       \centering
  \includegraphics[width=5in]{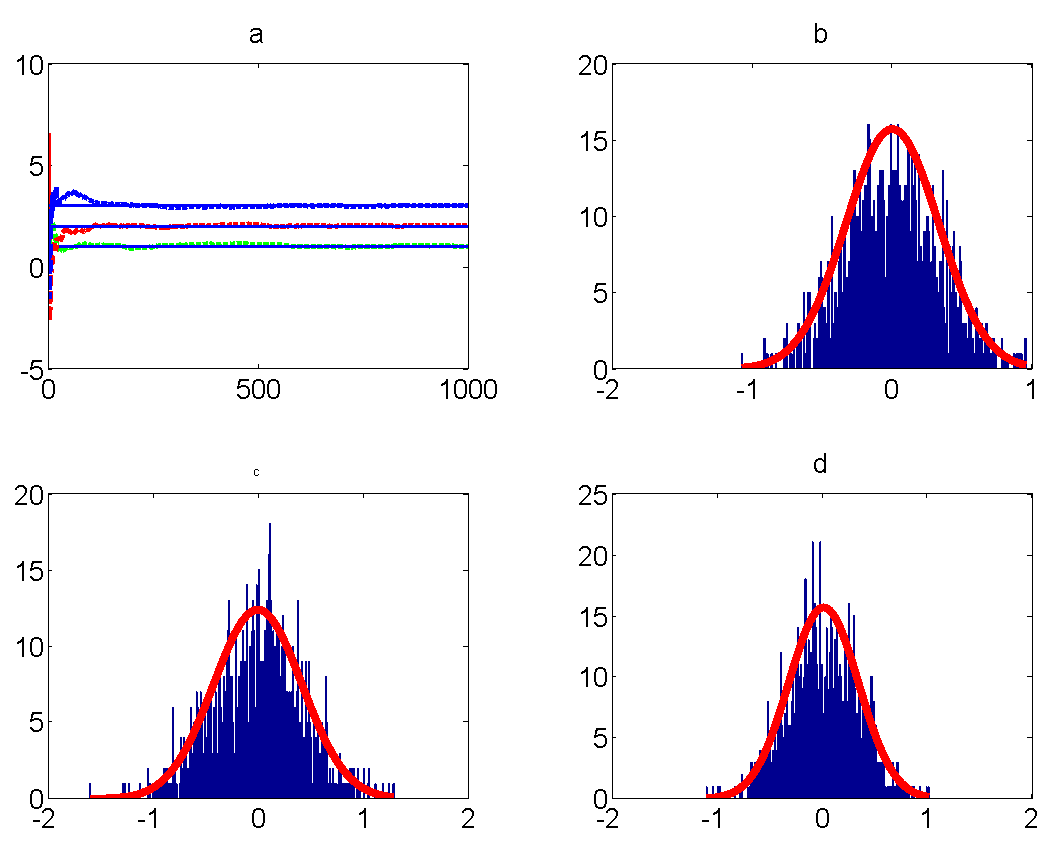}
   \caption{  Estimates $x_{2,k}$,   histographs and    limit distributions  for       $(x_{2,k}-x^*)/\sqrt{\gamma_k} $ at  $k=1000$} \label{GTwo}
\end{figure}

\begin{figure}
       \centering
  \includegraphics[width=5in]{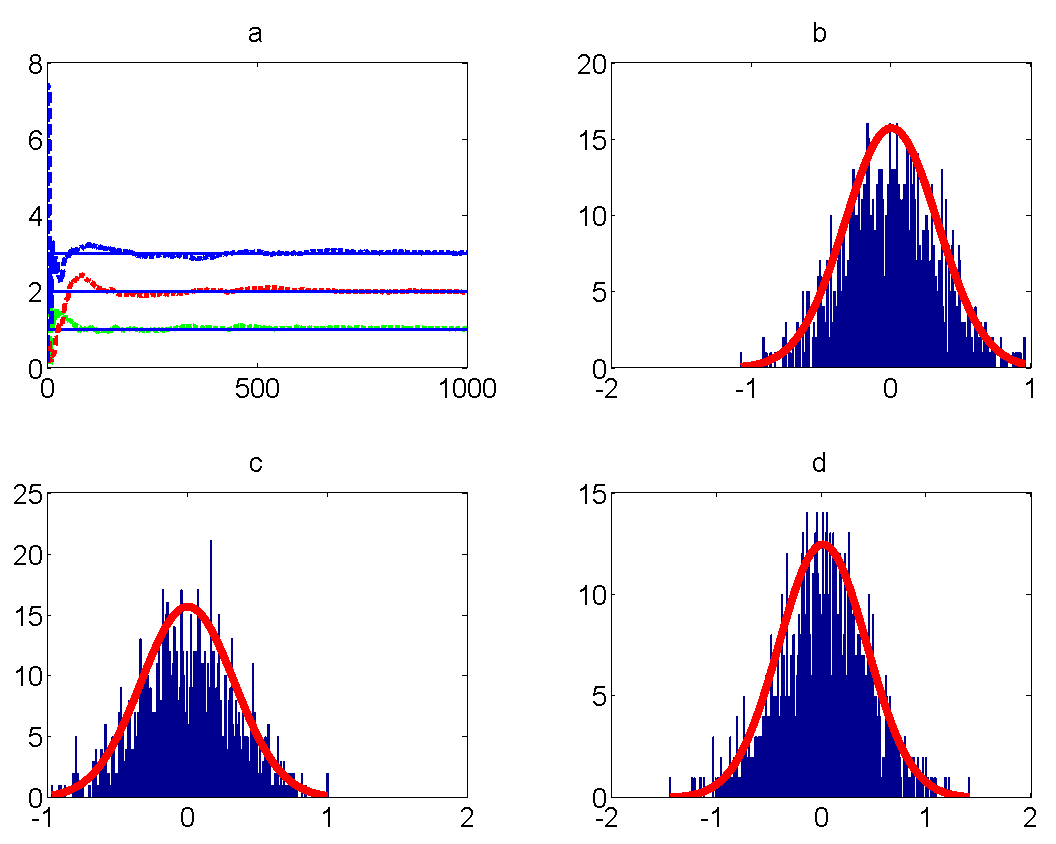}
    \caption{ Estimates $x_{3,k}$,   histographs and    limit distributions for      $(x_{3,k}-x^*)/\sqrt{\gamma_k} $ at $k=1000$} \label{GThree}
\end{figure}

 The algorithm \eqref{algorithm1} is calculated   for  $1000$  independent samples with  $k \leq 1000$.
   For $i=1,2,3,$ the estimates $x_{i,k}$ and the  histographs for  each component of 
   $(x_{i,k}-x^*)/\sqrt{\gamma_k}$ at time $1000$ are shown in Figs. \ref{GOne},   \ref{GTwo}, and   \ref{GThree}, respectively.   We use the normal distribution  to fit the  $1000$  independent samples
   for   each component of  $(x_{i,k}-x^*)/\sqrt{\gamma_k},~i=1,2,3$ with $k=1000$.  It  is shown  that the 
   data are fitted with  the  normal distribution by the Kolmogrov-Smirnov test  with  the significance level $\alpha=0.05.$    
   Fig. \ref{GOne}-a demonstrates estimates given by agent $1$ for components of $x^*=(1,2,3)$, where the real lines denote  true values, while the dashed lines are their estimates. The estimation errors  
   $(x_{i,k}-x^*)/\sqrt{\gamma_k}$   are presented in   Figs. \ref{GOne}-b,  \ref{GOne}-c,  \ref{GOne}-d, where the histographs are given by errors of 1000 samples at  time $k=1000$,  which are fitted by Gaussian densities. 
Figs.   \ref{GTwo} and     \ref{GThree} are for agents 2 and 3, respectively.

       \section{Conclusions}\label{sec:Conclusion}
       In this work, a   stochastic approximation based  distributed primal-dual algorithm is proposed   to  solve  the    distributed  constrained  stochastic   optimization  problem over random networks with imperfect communications.        
       The local estimates derived at all agents all shown to a.s. reach a consensus belonging to the optimal solution set.
  Besides,    we established conditions  for the unconstrained problem,  under which  the
       asymptotic normality and asymptotic efficiency of  the proposed  algorithm are established. 
       The  influence  on  the  convergence  rate   of the   
       conditional  covariance matrices   of  communication noises and gradient errors,   
       properties  of the cost function like   gradients and  Hessian  matrices  at the optimal point,  as well  as the  random graphs and its    mean graph  is demonstrated in the paper.

\appendices
\numberwithin{equation}{section}
\section{Some Results on Stochastic Approximation}

To ease reading, some results on non-neagetive super-martingales  \cite{supermds}   and  some information from stochastic approximation    \cite{Chen2002}  \cite{Kushner} are cited below.

\begin{lem} \cite{supermds} \label{martingale} 
Let $(\Omega,\mathcal{F},\mathbb{P})$ be a probability space and $\mathcal{F}_0\subset \mathcal{F}_1 \subset\cdots $
be a sequence of sub-$\sigma$-algebras of $\mathcal{F}$. Let $\{d_k\}$ and
$\{w_k\}$ be nonnegative  $\mathcal{F}_{k}$-measurable  random variables such that 
$$ E[d_{k+1}| \mathcal{F}_{k}] \leq (1+\alpha_k) d_k + w_k  ,$$
where $\alpha_k \geq 0$ and    $\sum_{k=1}^{\infty}\alpha_k <\infty$.
If  $\sum_{k=1}^{\infty} w_k < \infty~a.s.$, then  $\{d_k\}$ converges  a.s.  \end{lem}

We now  introduce the convergence results for  the constrained stochastic  approximation algorithm \cite{Kushner}. Consider the following  recursion 
\begin{equation}
\label{constrained}
\theta_{k+1}=P_{\Phi}(\theta_k+\gamma_kY_k),
\end{equation}
where $\Phi \in \mathbb{R}^m$ is a convex constraint set. 
 We list the conditions to be used.  
 
B1:  $\sup_k E[ \| Y_k\|^2]<\infty~~a.s.$ 

B2: There is a function $g(\cdot)$ such that
$$ E_k[Y_k]=E[Y_k|\theta_0,Y_i, i<k]=g(\theta_k)~~a.s.$$

B3: $g(\cdot)$ is continuous.

B4: $\theta_k$ is bounded  a.s.

\begin{lem}  \cite[Theorem 5.2.3] {Kushner} \label{CSA} 
Let $\{\theta_k\}$ be generated by \eqref{constrained}. 
   Assume  that the convex set  $\Phi$   satisfies  the same condition  as Assumption \ref{ass-set}-c   imposed on $\Omega_i$.  Let  B1-B4, and Assumption \ref{ass-stepsize} hold.  Then  $\theta_k$ converges a.s. to the limit  set of the following projected ODE   \cite{Kushner} in $\Phi$:
$$\dot{\theta}=g(\theta)+z ,~~z(t)\in -N_{\Phi} (\theta(t)),$$
where   $z(\cdot) $ is the projection or constraint  term,   the minimum force needed  to keep $\theta(\cdot)$ in $\Phi.$ 
\end{lem}

We   introduce asymptotic properties of the sequence $\{\theta_k\}$ generated by   the following recursion:
\begin{equation}\label{recursion}
\theta_{k+1}=\theta_k+\gamma_k g( \theta_k)+\gamma_ke_k.
\end{equation} 
We need the  following conditions.

 C0  There exists a continuously differentiable function $v(\cdot) $ such that
 $$g(x)^T \nabla v(x)<0~~\forall x \neq \mathbf{0}.$$
 
C1'  $  \theta_k$ is bounded a.s.

 C1   $\lim\limits_{k \rightarrow \infty} \theta_k=\mathbf{0}~~a.s.$

 C2' $\sum_{k=1}^{\infty} \gamma_k e_{k}<\infty~a.s.$
 
 C2  The noise  sequence  $\{e_{k}\}$ can be decomposed into two parts $e_{k}= \varepsilon_k+\nu_k$ such that 
 \begin{equation}\label{nuk}
\nu_k=o(\sqrt{\gamma_k})~~a.s.,
\end{equation} and   $\{\varepsilon_k,\mathcal{F}_k\}$ is an mds satisfying conditions: 
\begin{align}
&  E[\varepsilon_k| \mathcal{F}_{k-1}]=\mathbf{0},~~\sup_k E[\|\varepsilon_k\|^2| \mathcal{F}_{k-1}] \leq \sigma  
\textrm{~with }  \sigma \textrm{~being a constant}, \label{A-1}
\\& \lim_{k  \rightarrow \infty}  E[\varepsilon_k\varepsilon_k^T| \mathcal{F}_{k-1}]= \lim_{k  \rightarrow \infty}  E[\varepsilon_k\varepsilon_k^T]=S_0~~a.s.,\label{A-2}
\\& \lim_{a \rightarrow \infty} \sup_k E[\|\varepsilon_k\|^2I_{[\|\varepsilon_k\|>a]}]=0. \label{A-3}
\end{align}

 C3' $g(\cdot)$ is measurable and locally bounded. 
 
 C3 $g(\cdot)$ is measurable and locally bounded. As $\theta \rightarrow \mathbf{0},$
 $$\| g(\theta)-F\theta\|\leq c\|\theta\|^2,$$
 where $c>0$ and $F$ is stable. 
 
   \begin{lem} \label{Stable0}\cite[Theorem 2.2.1]{Chen2002} Let  $\{\theta_k\}$ by generated by  \eqref{recursion} with an arbitrary initial value $\theta_0$.     Let  Assumption \ref{ass-stepsize1},  and C0, C1', C2', and C3' hold.   Then
 $$\lim\limits_{k \rightarrow \infty} \theta_k=\mathbf{0}~~a.s.$$
\end{lem}

 \begin{lem} \label{Stable1} Let  $\{\theta_k\}$ by generated by  \eqref{recursion}.     Let  Assumption \ref{ass-stepsize1},  and C0, C1, C2, and C3 hold.   Then

   i)  $\frac{1}{\sqrt{\gamma_k}} \theta_k$ is asymptotically normal:
 $$\frac{1}{\sqrt{\gamma_k}} \theta_k \xlongrightarrow [k  \rightarrow \infty]  {d}N(\textbf{0},S),$$
where $S=\int_{0}^{\infty} e^{Ft}S_0 e^{F^Tt} dt$; 
 
 ii)     $\bar{\theta}_k$ is asymptotically efficient: 
 $$ \sqrt{k} \bar{\theta}_k \xlongrightarrow [k  \rightarrow \infty]  {d}N(\textbf{0},S),$$
where $S=F^{-1}S^0 (F^{-1})^T,$ and  $\bar{\theta}_k=\frac{1}{k}\sum_{p=1}^k \theta_p $.
\end{lem} 
 \begin{rem}
 Lemma \ref{Stable1} i) is  \cite[Theorem 3.3.2]{Chen2002}  for the   case:   $r=0, \beta=1, \alpha=0, x^0=\mathbf{0}$.  Since 
 the noise   sequence  $\{e_{k}\}$ satisfies C2, by \cite[Remarks 3.4.1 and  3.4.2]{Chen2002}   it is seen that  A3.4.3 in \cite{Chen2002}  holds. Then by \cite[Theorem 3.4.2]{Chen2002}   with  $  \beta=1,~x^0=\mathbf{0}$  the assertion of 
 Lemma \ref{Stable1} ii) follows.  
 \end{rem} 
  
   \end{document}